\documentclass[a4paper, 12pt]{article}
\usepackage[top=1in, bottom=1in, left=1.25in, right=1.25in]{geometry}
\usepackage[dvips]{graphicx}
\usepackage{natbib}
\usepackage{authblk}
\usepackage{graphicx}
\usepackage{amssymb,amsmath,amsthm,latexsym}
\usepackage{ascmac}
\usepackage{mathtools}
\usepackage{booktabs}
\usepackage{threeparttable}

\DeclarePairedDelimiterX\Set[2]{\{}{\}}{\mspace{2mu}{#1}\;\delimsize|\;{#2}\mspace{2mu}}
\mathtoolsset{showonlyrefs}
\usepackage{amsfonts}
\usepackage{amssymb}
\usepackage{bm}
\usepackage{bbm}
\usepackage{pifont}
\usepackage{fancybox}
\usepackage{enumitem}
\usepackage{color}
\usepackage{arydshln}
\usepackage[T1]{fontenc}
\usepackage{lscape}
\usepackage{caption}

\usepackage[, 
	colorlinks = true, 
	anchorcolor=black,
	linkcolor = blue,
	citecolor = blue, 
	filecolor  = blue, 
	urlcolor   = blue
]{hyperref}
\usepackage{mathrsfs}
\usepackage{comment}
\usepackage[titletoc]{appendix}
\usepackage{hyperref}

\usepackage[utf8]{inputenc} 
\usepackage{float}

\theoremstyle{plain}
\newtheorem{thm}{Theorem}
\newtheorem{lem}{Lemma}

\theoremstyle{definition}

\theoremstyle{remark}

\renewcommand{\vec}[1]{\mbox{\boldmath ${#1}$}}

\allowdisplaybreaks

\title{\bf Edgeworth expansions in sphericity testing under high-dimensional two-step monotone incomplete data and their computable error bounds}
\author[1]{Tetsuya Sato}
\author[2,3]{Tomoyuki Nakagawa}

\affil[1]{Department of Applied Mathematics, Graduate School of Science, Tokyo University of Science\footnote{1-3 Kagurazaka, Shinjuku-ku Tokyo 162-8601, Japan}}
\affil[2]{School of Data Science, Meisei University\footnote{2-1-1 Hodokubo, Hino, Tokyo 191-8506, Japan}}
\affil[3]{Statistical Mathematics Unit, RIKEN Center for Brain Science\footnote{2-1 Hirosawa Wako City, Saitama 351-0198 Japan}}
\date{Last update: \today}

\begin{document}
\maketitle
\begin{abstract}
In this paper, we consider the sphericity test for a one-sample problem under high-dimensional two-step monotone incomplete data. 
Existing asymptotic expansions for the null distributions of the likelihood ratio test (LRT) statistic and modified LRT statistic deteriorate in high-dimensional settings. 
Therefore, we derive an Edgeworth expansion for the null distribution of the LRT statistic in such settings and obtain its computable error bounds. 
Furthermore, we demonstrate that our proposed Edgeworth expansion provide better approximation accuracy than the existing asymptotic expansion. 
We also conduct numerical experiments using Monte Carlo simulations to evaluate the maximum absolute error between the distribution function of the standardized test statistic and the Edgeworth expansion for the null distribution of the LRT statistic, as well as to assess the performance of the computable error bounds. Finally, by deriving a Cornish-Fisher expansion for the null distribution of the LRT statistic in high-dimensional settings, we evaluate the empirical Type I error rates via Monte Carlo simulations.
\end{abstract}


\medskip


\baselineskip21.5pt
\section{Introduction}
\label{sec:intro}
\begin{sloppypar}
This paper studies the sphericity test for a one-sample problem under high-dimensional two-step monotone incomplete data.  
The sphericity hypothesis is given by
\begin{align} 
H_0 : \bm{\Sigma} = \sigma^2 \bm I_p \mbox{ vs.}\  H_1 : \bm{\Sigma} \neq \sigma^2 \bm I_p,\label{eqn:a1}
\end{align}
where $\bm{\Sigma}$ denotes the covariance matrix and $\sigma^2$ is an unknown positive parameter.
The null hypothesis $H_0$ corresponds to an isotropic covariance structure, in which all variables have equal variance and are mutually uncorrelated.
Testing sphericity plays a fundamental role in multivariate analysis, as rejection of $H_0$ signals meaningful covariance patterns and motivates the use of methods such as principal component analysis.
\end{sloppypar} 
For the sphericity testing problem under a multivariate normal distribution, \cite{mauchly1940significance} derived the likelihood ratio (LR) and its null moments under complete data. \cite{mauchly1940significance} also noted that the asymptotic null distribution of the likelihood ratio test (LRT) statistic follows a chi-square distribution as shown by \cite{Wilks1938}. Building on this work, \cite{gleser1966note} proved that Mauchly's test is unbiased and derived the asymptotic distribution of the test statistic under $H_1$. 
Subsequently, \cite{john1971} derived the locally most powerful invariant test for sphericity, demonstrating that a statistic based on the sample covariance's eigenstructure is optimal under these assumptions. 
Following this, \cite{john1972distribution} obtained the exact null distribution of a sphericity statistic and provided its asymptotic and approximate null distributions. 
In addition, \cite{muirhead1982aspects} and \cite{Anderson2003aspects} provided the asymptotic expansion for the null distribution of the modified LRT statistic under complete data. These works typically assume a low-dimensional framework and rely on large-sample asymptotics. 
More recently, in high-dimensional settings where the dimension exceeds the sample size, \cite{ledoit2002some} studied a sphericity statistic shown in \cite{john1971} rather than the LRT statistic and \cite{Srivastava2006} derived the asymptotic expansion for the null distribution of the modified LRT statistic as the dimension tends to infinity and the sample size is fixed. In addition, \cite{wang2013sphericity} proposed a corrected LRT statistic rewritten in terms of eigenvalues, which differs from the LRT statistic, as well as a corrected John's test statistic, and derived normal approximations for their null distributions under high-dimensional asymptotics.

\begin{sloppypar}
Edgeworth expansions and their computable error bounds in high-dimensional settings have been derived in several contexts. These include the work of \cite{wakaki2006edgeworth} on Wilks' lambda statistic introduced by \cite{Wilks1932}; \cite{wakaki2007error} on the null distributions of three tests including the sphericity test and the non-null distribution of one of these tests; \cite{akita2010high} on a test statistic for independence; \cite{kato2010asymptotic} on testing the intraclass correlation structure; and \cite{wakaki2018computable} on an LRT for additional information hypothesis in canonical correlation analysis.
Furthermore, \cite{Ulyanov2016} established computable error bounds for generalized Cornish-Fisher expansions of quantiles. These bounds were derived assuming that error bounds for the corresponding Edgeworth expansions are available. Additionally, \cite{fujikoshi2020non} introduced large-sample asymptotic expansions and error bounds for general high-dimensional asymptotic theory.
The derivation of such error bounds is crucial, as they allow for assessing the accuracy of asymptotic approximations and validating their use for finite samples, which is particularly vital in high-dimensional settings.
However, studies on Edgeworth expansions and their computable error bounds for incomplete data are scarce, with most existing work focusing on the complete-data case.
\end{sloppypar}

\begin{sloppypar}
In contrast, in monotone incomplete-data settings, \cite{batsidis2006k-step} derived the LR for general $k$-step monotone incomplete data under an elliptical distribution and \cite{chang2010finite} discussed the null moments of the LR under two-step monotone incomplete data. 
\cite{Sato2025sphericity} derived asymptotic expansions for the null distributions of the LRT statistic and modified LRT statistic under two-step and $k$-step monotone incomplete data for the large sample size using the general distribution proposed by \cite{box1949general}. 
Furthermore, \cite{satomulti} extended \cite{Sato2025sphericity} to a multi-sample problem. However, to the best of our knowledge, Edgeworth expansions and their error bounds for monotone incomplete data in high-dimensional settings have not been addressed in the literature. 
Nevertheless, incomplete data frequently arise in real data analysis, making it practically important to account for incompleteness. 
Therefore, we derive our Edgeworth expansion for the null distribution in sphericity testing under two-step monotone incomplete data, together with the computable error bounds to assess its validity.
\end{sloppypar}

The remainder of this paper is organized as follows. 
In Section \ref{sec:LR}, we describe the LR for the sphericity test \eqref{eqn:a1} under a multivariate normal distribution and asymptotic expansions for the null distributions of the LRT statistic and modified LRT statistic for the large sample size under two-step monotone incomplete data. 
In Section \ref{sec:Edgeworth}, we obtain the Edgeworth expansion for the null distribution of the LRT statistic by deriving the characteristic function. 
In Section \ref{sec:Error bounds}, we derive a uniform bound for the error of the Edgeworth expansion and some upper bounds. 
In Section \ref{subsec:Type I error rates}, we use Monte Carlo simulations to compare the approximation accuracy of the Edgeworth expansion and the large-sample asymptotic expansion for the null distribution of the LRT statistic. In Section \ref{subsec:Example}, we numerically evaluate the maximum absolute error (MAE) between the distribution function of the standardized test statistic and the Edgeworth expansion for the null distribution using Monte Carlo simulations and perform numerical experiments to evaluate the error bounds. In Section \ref{subsec:Cornish-Fisher}, we derive a high-dimensional Cornish-Fisher expansion for the null distribution of the LRT statistic and evaluate the empirical Type I error rates and power through Monte Carlo simulations. Finally, in Section \ref{sec:Conclusions}, the conclusions are presented. The proofs are provided in \ref{app:A}--\ref{app:C}, while the additional simulation results are presented in \ref{app:D}.
\section{Likelihood ratio test for sphericity}
\label{sec:LR}
\begin{sloppypar}
Let $\bm{x}_{1}, \ldots, \bm{x}_{N_1}$ be independently and identically distributed (i.i.d.) multivariate normal $N_{p}(\bm{\mu}, \bm{\Sigma})$ and let $\bm{x}_{1N_1+1},\ldots,\bm{x}_{1N}$ be i.i.d. multivariate normal $N_{p_{1}}(\bm{\mu}_{1}, \bm{\Sigma}_{11})$, where
\end{sloppypar}
\begin{align*}
\bm\mu=\begin{pmatrix}
\bm\mu_1\\
\bm\mu_2
\end{pmatrix},
\bm\Sigma=\begin{pmatrix}
\bm\Sigma_{11} &\bm\Sigma_{12}\\
\bm\Sigma_{21} &\bm\Sigma_{22}
\end{pmatrix}.
\end{align*}
Here, $\bm\mu_i$ is a $p_i\times 1$ vector, $\bm\Sigma_{ij}$ is a $p_i\times p_j$ matrix for $i,j=1,2$, and $N=N_1+N_2$, $p=p_1+p_2$. We assume that $N_1>p$. 
Then, the two-step monotone incomplete data structure can be expressed as
\begin{align*}
\left( \ 
\begin{array}{|c|c|c|c|c|c|}
\cline{1-6}
\multicolumn{5}{|c}{\rule{0pt}{2.2ex}\bm{x}_{11}^\top}&
\multicolumn{1}{c|}{\rule{0pt}{2.2ex}\bm{x}_{21}^\top}\\
\multicolumn{5}{|c}{\vdots}&
\multicolumn{1}{c|}{\vdots}\\
\multicolumn{5}{|c}{\bm{x}_{1N_1}^\top}&
\multicolumn{1}{c|}{\bm{x}_{2N_1}^\top}\\
\cline{6-6}
\multicolumn{5}{|c|}{\bm{x}_{1N_1+1}^\top}&
\multicolumn{1}{c}{\ast \hspace{0.25cm} \cdots \hspace{0.25cm} \ast}\\
\multicolumn{5}{|c|}{\vdots}&
\multicolumn{1}{c}{\vdots \hspace{1.15cm} \vdots}\\
\multicolumn{5}{|c|}{\bm{x}_{1N}^\top}&
\multicolumn{1}{c}{\ast \hspace{0.25cm} \cdots \hspace{0.25cm} \ast}\\[1pt]
\cline{1-5}
\end{array}
\ \right),
\end{align*}
where $\vec{x}_{j}=\begin{pmatrix}
\bm{x}_{1j}^\top,\bm{x}_{2j}^\top
\end{pmatrix}^\top,j=1,\ldots,N_1$ and `$*$' indicates unobserved data.
As can be seen from this notation, $\bm{x}_{11}, \ldots, \bm{x}_{1N_1}$ and $\bm{x}_{1N_1+1}, \ldots, \bm{x}_{1N}$ share common parameters $\bm{\mu}_1$ and $\bm{\Sigma}_{11}$, and $\bm{x}_{1N_1+1}, \ldots, \bm{x}_{1N}$ contains information about them. 
The LR for the sphericity test \eqref{eqn:a1} under two-step monotone incomplete data is given by
\begin{align*}
\lambda=|\widehat{\bm\Sigma}_{11}|^{\tfrac{N}{2}}|\widehat{\bm\Sigma}_{22\cdot1}|^{\tfrac{N_1}{2}}(\widetilde{\sigma}^2)^{-\tfrac{Np_1+N_1p_2}{2}},
\end{align*}
where $\bm\Sigma_{22\cdot1}=\bm\Sigma_{22}-\bm\Sigma_{21}\bm\Sigma_{11}^{-1}\bm\Sigma_{12}$ and 
$\widehat{\bm\Sigma}=\begin{pmatrix}
\widehat{\bm\Sigma}_{11} &\widehat{\bm\Sigma}_{12}\\
\widehat{\bm\Sigma}_{21} &\widehat{\bm\Sigma}_{22}
\end{pmatrix}$ 
is the maximum likelihood estimator (MLE) of $\bm\Sigma$, as derived by \cite{Anderson1985mle} and \cite{kanda1998some}, $\widetilde{\sigma}^2$ is the MLE of $\sigma^2$ under $H_0$ given by \cite{chang2010finite}, and
\begin{align*}
\widehat{\bm\Sigma}_{11}&=\dfrac{1}{N}\left(\bm W_{11}^{(1)}+\bm W_{11}^{(2)}\right),
~\widehat{\bm\Sigma}_{12}=\widehat{\bm\Sigma}_{21}^\top=\widehat{\bm\Sigma}_{11}(\bm W_{11}^{(1)})^{-1}\bm W_{12}^{(1)},\\
\widehat{\bm\Sigma}_{22}&=\dfrac{1}{N_1}\bm W_{22\cdot1}^{(1)}+\widehat{\bm\Sigma}_{21}\widehat{\bm\Sigma}_{11}^{-1}\widehat{\bm\Sigma}_{12},
~\widehat{\bm\Sigma}_{22\cdot1}=\widehat{\bm\Sigma}_{22}-\widehat{\bm\Sigma}_{21}\widehat{\bm\Sigma}_{11}^{-1}\widehat{\bm\Sigma}_{12},\\
\widetilde{\sigma}^2&=\dfrac{1}{Np_1+N_1p_2}\left(\mathrm{tr}(\bm W_{11}^{(1)}+\bm W_{11}^{(2)})+\mathrm{tr}\bm W_{22}^{(1)}\right),\\
\overline{\bm{x}}_1^{(1)}&=\dfrac{1}{N_1}\sum_{j=1}^{N_1}\bm{x}_{1j},
~\overline{\bm{x}}_2^{(1)}=\dfrac{1}{N_1}\sum_{j=1}^{N_1}\bm{x}_{2j},
~\overline{\bm{x}}_1^{(2)}=\dfrac{1}{N_2}\sum_{j=N_1+1}^{N}\bm{x}_{1j},\\
\bm{W}_{11}^{(1)}&=\sum_{j=1}^{N_1} (\bm{x}_{1j}-\overline{\bm{x}}_{1}^{(1)})(\bm{x}_{1j}-\overline{\bm{x}}_{1}^{(1)})^\top,
~\bm{W}_{22}^{(1)}=\sum_{j=1}^{N_1} (\bm{x}_{2j}-\overline{\bm{x}}_{2}^{(1)})(\bm{x}_{2j}-\overline{\bm{x}}_{2}^{(1)})^\top,\\  
\bm{W}_{12}^{(1)}&=(\bm{W}_{21}^{(1)})^\top=\sum_{j=1}^{N_1} (\bm{x}_{1j}-\overline{\bm{x}}_{1}^{(1)})(\bm{x}_{2j}-\overline{\bm{x}}_{2}^{(1)})^\top,\\
\bm{W}_{11}^{(2)}&=\sum_{j=N_1+1}^{N} (\bm{x}_{1j}-\overline{\bm{x}}_{1}^{(2)})(\bm{x}_{1j}-\overline{\bm{x}}_{1}^{(2)})^\top+\dfrac{N_1N_2}{N}(\overline{\bm{x}}_{1}^{(1)}-\overline{\bm{x}}_{1}^{(2)})(\overline{\bm{x}}_{1}^{(1)}-\overline{\bm{x}}_{1}^{(2)})^\top,\\
\bm{W}_{22\cdot 1}^{(1)}&=\bm{W}_{22}^{(1)}-\bm{W}_{21}^{(1)}(\bm{W}_{11}^{(1)})^{-1}\bm{W}_{12}^{(1)}.
\end{align*}
Furthermore, following \cite{chang2010finite}, the LR can be written as
\begin{align*}
\lambda=\dfrac{\left|\dfrac{1}{N}\left(\bm W_{11}^{(1)}+\bm W_{11}^{(2)}\right)\right|^{\tfrac{N}{2}}\left|\dfrac{1}{N_1}\bm W_{22\cdot1}^{(1)}\right|^{\tfrac{N_1}{2}}}{\biggl\{\dfrac{1}{Np_1+N_1p_2}\left(\mathrm{tr}(\bm W_{11}^{(1)}+\bm W_{11}^{(2)})+\mathrm{tr}\bm W_{22}^{(1)}\right)\biggr\}^{\tfrac{Np_1+N_1p_2}{2}}}.
\end{align*}
Based on this representation, the asymptotic expansions for the null distributions of the LRT statistic and modified LRT statistic are given in the following theorem.
\begin{thm}[\cite{Sato2025sphericity}]\label{thm:1}
Adapting the result of \cite{Sato2025sphericity}, when $\tau_1=N_1/N\to\delta_1\in(0,1]~(N_1\to\infty)$, the asymptotic expansions for the null distributions of the LRT statistic and modified LRT statistic are given for large $N_1$ as
\begin{align}
\begin{split}\label{eqn:a2}
\mathrm{Pr}(-2\log\lambda\leq x)=&~G_{f}(x)+\dfrac{\beta}{N}\bigl\{G_{f+2}(x)-G_{f}(x)\bigr\}\\
&+\dfrac{\gamma}{N^2}\bigl\{G_{f+4}(x)-G_{f}(x)\bigr\}+O(N^{-3}),
\end{split}\\
\mathrm{Pr}(-2\rho\log\lambda\leq x)=&~G_{f}(x)+\dfrac{\gamma^{*}}{M^2}\bigl\{G_{f+4}(x)-G_{f}(x)\bigr\}+O(M^{-3}),\label{eqn:a3}
\end{align}
where $G_k(x)$ is the distribution function of the $\chi^{2}$ distribution with $k$ degrees of freedom, $f=(p+2)(p-1)/2$, $M=\rho N$, and
\begin{align*}
\beta=&~\dfrac{1}{24}\biggl[p_1(2p_1^2+9p_1+11)\\
&+\dfrac{1}{\tau_1}\bigl\{p_2(2p_2^2+9p_2+11)+6p_1p_2(p+3)\bigr\}-\dfrac{2}{p_1+\tau_1p_2}(3p^2+6p+2)\biggr],\\
\gamma=&~\dfrac{1}{48}\Biggl[p_1(p_1+1)(p_1+2)(p_1+3)+\dfrac{1}{\tau_1^2}\biggl(p_2(p_2+1)(p_2+2)(p_2+3)\\
&+2p_1p_2\bigl\{(p_2+1)(2p+p_1+7)+2(p_1+1)(p_1+2)\bigr\}\biggr)\\
&-\dfrac{4}{(p_1+\tau_1p_2)^2}p(p+1)(p+2)\Biggr],\\
\rho=&~1-\dfrac{4}{(p+2)(p-1)N}\beta,\\
\gamma^{*}=&-\dfrac{2}{(p+2)(p-1)}\beta^2+\gamma.
\end{align*}
\end{thm}
Theorem \ref{thm:1} provides the asymptotic results for fixed $p$. However, the accuracy of \eqref{eqn:a2} and \eqref{eqn:a3} is observed to deteriorate as $p$ increases, as shown by Monte Carlo simulations in Section \ref{subsec:Cornish-Fisher} and \ref{app:D}. Therefore, we derive an Edgeworth expansion for the null distribution under high-dimensional two-step monotone incomplete data, with the complete data as a special case, together with its corresponding error bounds, where $n=N-1,n_1=N_1-1,n\ge n_1>p>p_1$, and $0<\tau_1\leq 1$.
\section{Edgeworth expansions}
\label{sec:Edgeworth}
We derive our Edgeworth expansion for the null distribution of the LRT statistic. In this section, when reduced to the complete data case where $N_2=0$, our results agree with those reported in \cite{wakaki2007error} for the high-dimensional sphericity test.

For $h=0,1,2,\ldots$, following \cite{chang2010finite}, the $h$-th null moment of $\lambda$ is expressed as
\begin{align*}
\mathrm{E}[\lambda^{h}]=&\left[\dfrac{(Np_1+N_1p_2)^{\tfrac{Np_1+N_1p_2}{2}}}{N^{\tfrac{Np_1}{2}}N_1^{\tfrac{N_1p_2}{2}}}\right]^h\\
&\times\displaystyle\prod_{\ell=1}^{p_1}\dfrac{\Gamma\left[\dfrac{1}{2}(n-p_1+\ell)+\dfrac{1}{2}Nh\right]}{\Gamma\left[\dfrac{1}{2}(n-p_1+\ell)\right]}\prod_{\ell=1}^{p_2}\dfrac{\Gamma\left[\dfrac{1}{2}(n_1-p+\ell)+\dfrac{1}{2}N_1h\right]}{\Gamma\left[\dfrac{1}{2}(n_1-p+\ell)\right]}\\
&\times\dfrac{\Gamma\left[\dfrac{1}{2}(np_1+n_1p_2)\right]}{\Gamma\left[\dfrac{1}{2}(np_1+n_1p_2)+\dfrac{1}{2}(Np_1+N_1p_2)h\right]}.
\end{align*}
Using this formula, the characteristic function of $-(2/N)\log\lambda$ can be written as 
\begin{align*}
\varphi_{\lambda}(t)
=&~\mathrm{E}\bigl[\lambda^{-\tfrac{2it}{N}}\bigr]\\
=&~\Biggl(\dfrac{(p_1+\tau_1p_2)^{p_1+\tau_1p_2}}{\tau_1^{\tau_1p_2}}\Biggr)^{-it}\\
&\times\prod_{\ell=1}^{p_1}\dfrac{\Gamma\left[\dfrac{1}{2}(n-p_1+\ell)-it\right]}{\Gamma\left[\dfrac{1}{2}(n-p_1+\ell)\right]}\prod_{\ell=1}^{p_2}\dfrac{\Gamma\left[\dfrac{1}{2}(n_1-p+\ell)-\tau_1it\right]}{\Gamma\left[\dfrac{1}{2}(n_1-p+\ell)\right]}\\
&\times\dfrac{\Gamma\left[\dfrac{1}{2}(np_1+n_1p_2)\right]}{\Gamma\left[\dfrac{1}{2}(np_1+n_1p_2)-(p_1+\tau_1p_2)it\right]}.
\end{align*}
For $s\ge2$, the cumulants of $-(2/N)\log\lambda$ are given by
\begin{align*}
\kappa_{\lambda}^{(1)}=&-(p_1+\tau_1p_2)\log(p_1+\tau_1p_2)+\tau_1p_2\log\tau_1
+(p_1+\tau_1p_2)\psi^{(0)}\biggl(\dfrac{1}{2}(np_1+n_1p_2)\biggr)\\
&-\displaystyle\sum_{\ell=1}^{p_1}\psi^{(0)}\biggl(\dfrac{1}{2}(n-p_1+\ell)\biggr)-\tau_1\displaystyle\sum_{\ell=1}^{p_2}\psi^{(0)}\biggl(\dfrac{1}{2}(n_1-p+\ell)\biggr),\\
\kappa_{\lambda}^{(s)}=&~(-1)^{s-1}\Biggl[(p_1+\tau_1p_2)^{s}\psi^{(s-1)}\biggl(\dfrac{1}{2}(np_1+n_1p_2)\biggr)\\
&-\displaystyle\sum_{\ell=1}^{p_1}\psi^{(s-1)}\biggl(\dfrac{1}{2}(n-p_1+\ell)\biggr)-\tau_1^s\displaystyle\sum_{\ell=1}^{p_2}\psi^{(s-1)}\biggl(\dfrac{1}{2}(n_1-p+\ell)\biggr)\Biggr].
\end{align*}
Here, $\psi^{(s)}(a)$ denotes the polygamma function. Equivalently,
\begin{align*}
\psi^{(s)}(a)=\left.\frac{d^{s+1}}{dz^{s+1}}\log\Gamma[z]\right|_{z=a}
=\begin{cases}
-C+\displaystyle\sum_{k=0}^{\infty}\biggl(\dfrac{1}{1+k}-\dfrac{1}{k+a}\biggr) & \text{($s=0$)}, \\
    \displaystyle\sum_{k=0}^{\infty} \dfrac{(-1)^{s+1}s!}{(k+a)^{s+1}} & \text{($s\ge1$),}
  \end{cases}
\end{align*}
where $C$ is the Euler constant. Let the standardized test statistic be
\begin{align}
T=\dfrac{-\tfrac{2}{N}\log\lambda-\kappa_{\lambda}^{(1)}}{(\kappa_{\lambda}^{(2)})^{\tfrac{1}{2}}},\label{eqn:a3.5}
\end{align}
and denote its $s$-th standardized cumulant by
\begin{align*}
\widetilde{\kappa}_{T}^{(s)}=\dfrac{\kappa_{\lambda}^{(s)}}{(\kappa_{\lambda}^{(2)})^{\tfrac{s}{2}}},
\end{align*}
for $s\ge3$. 
Then, the upper bounds for the standardized cumulants are given by the following lemma.
\begin{lem}\label{lem:1}
For $s\ge0$, define $m$ and $b_s$ by
\begin{align*}
m=&~\dfrac{n_1-p-\tfrac{1}{2}}{2}(\kappa_{\lambda}^{(2)})^{\tfrac{1}{2}},\\
b_s=&~\dfrac{2}{\kappa_{\lambda}^{(2)}(s+1)(s+2)(s+3)}\Biggl\{\biggl(\dfrac{n_1-p-\tfrac{1}{2}}{n-p_1-\tfrac{1}{2}}\biggr)^{s+1}-\biggl(\dfrac{n_1-p-\tfrac{1}{2}}{n-\tfrac{1}{2}}\biggr)^{s+1}\Biggr\}\\
&+\dfrac{2\tau_1^2}{\kappa_{\lambda}^{(2)}(s+1)(s+2)(s+3)}\Biggl\{\tau_1^{s+1}-\biggl(\tau_1\dfrac{n_1-p-\tfrac{1}{2}}{n_1-p_1-\tfrac{1}{2}}\biggr)^{s+1}\Biggr\}\\
&-\Biggl\{\dfrac{2}{\kappa_{\lambda}^{(2)}(s+3)n^2}+\dfrac{2(p_1+\tau_1p_2)}{\kappa_{\lambda}^{(2)}(s+2)(s+3)n}\Biggr\}\biggl(\dfrac{n_1-p-\tfrac{1}{2}}{n}\biggr)^{s+1}.
\end{align*}
Then it holds that
\begin{align}
0<\dfrac{\widetilde{\kappa}_{T}^{(s)}}{s!}<m^{-(s-2)}b_{s-3},\label{eqn:a4}
\end{align}
for $s\ge3$.
\end{lem}
See \ref{app:A} for the proof of Lemma \ref{lem:1}. We can check that $b_s$ is bounded and $m$ becomes large if $n_1$ and $p$ become large. The characteristic function of $T$ can be expanded as
\begin{align*}
\varphi_{T}(t)
=&\exp\Bigl(-\dfrac{t^2}{2}+\sum_{s=3}^{\infty}\dfrac{\widetilde{\kappa}_{T}^{(s)}}{s!}(it)^s\Bigr)\\
=&\exp\Bigl(-\dfrac{t^2}{2}\Bigr)\Biggl\{1+\sum_{k=1}^{\infty}(it)^{3k}\Bigl(\sum_{s=0}^{\infty}\dfrac{\widetilde{\kappa}_{T}^{(s+3)}}{(s+3)!}(it)^s\Bigr)^k\Biggr\}\\
=&\exp\Bigl(-\dfrac{t^2}{2}\Bigr)\Biggl\{1+\sum_{k=1}^{\infty}\dfrac{1}{k!}(it)^{3k}\sum_{j=0}^{\infty}\gamma_{k,j}(it)^j\Biggr\},
\end{align*}
where
\begin{align*}
\gamma_{k,j}=\sum_{s_1+\cdots+s_k=j}\dfrac{\widetilde{\kappa}_{T}^{(s_1+3)}\cdots\widetilde{\kappa}_{T}^{(s_k+3)}}{(s_1+3)!\cdots(s_k+3)!}.
\end{align*}
The inequality \eqref{eqn:a4} implies $\gamma_{k,j}=O(m^{-(j+k)})$ (see \ref{app:B} for the proof). Hence, we have
\begin{align*}
\varphi_{T}(t)=\varphi_{T,s}(t)+O(m^{-(s+1)})\qquad(m\to\infty),
\end{align*}
where
\begin{align}
\varphi_{T,s}(t)=\exp\Bigl(-\dfrac{t^2}{2}\Bigr)\Biggl\{1+\sum_{k=1}^{s}\dfrac{1}{k!}(it)^{3k}\sum_{j=0}^{s-k}\gamma_{k,j}(it)^j\Biggr\}.\label{eqn:a5}
\end{align}
Inverting \eqref{eqn:a5} yields the Edgeworth expansion for the null distribution of the LRT statistic up to the order $O(m^{-s})$:
\begin{align}
Q_s(x)=\Phi(x)-\phi(x)\Biggl\{\sum_{k=1}^{s}\dfrac{1}{k!}\sum_{j=0}^{s-k}\gamma_{k,j}h_{3k+j-1}(x)\Biggr\},\label{eqn:a6}
\end{align}
where $\Phi(x)$ and $\phi(x)$ are the distribution function and probability density function of the standard normal distribution, respectively, and $h_r(x)$ is the $r$-th order Hermite polynomial defined by
\begin{align*}
\phi(x)h_r(x)=(-1)^r\frac{d^{r}}{dx^{r}}\phi(x).
\end{align*}
In the next section, we derive the computable error bounds to assess the validity of the Edgeworth expansion $Q_s(x)$.
\section{Error bounds}
\label{sec:Error bounds}
Error bounds are essential to quantify the accuracy of the asymptotic approximation and to ensure its validity in finite samples, particularly in high-dimensional settings.
Using the inverse Fourier transformation, we obtain a uniform bound on the error of the Edgeworth expansion $Q_s(x)$.
\begin{align*}
\sup_{x}\bigl|\Pr(T\leq x)-Q_s(x)\bigr|&\leq\dfrac{1}{2\pi}\int_{-\infty}^{\infty}\dfrac{1}{|t|}\left|\varphi_{T}(t)-\varphi_{T,s}(t)\right|dt\\
&\leq\dfrac{1}{2\pi}\Bigl(I_1(v)+I_2(v)+I_3(v)\Bigr),
\end{align*}
where
\begin{align*}
I_1(v)&=\int_{-mv}^{mv}\dfrac{1}{|t|}\left|\varphi_{T}(t)-\varphi_{T,s}(t)\right|dt,\\
I_2(v)&=\int_{|t|>mv}\dfrac{1}{|t|}\left|\varphi_{T,s}(t)\right|dt,
I_3(v)=\int_{|t|>mv}\dfrac{1}{|t|}\left|\varphi_{T}(t)\right|dt,
\end{align*}
for some constant $v\in(0,1)$. We derive bounds for $I_1(v)$, $I_2(v)$, and $I_3(v)$.

First, we derive an upper bound $U_1(v)$ for $I_1(v)$. Let $L_1(v)$, $L_2(v)$, and $L_3(v)$ be
\begin{align*}
L_1(v)=&~\sum_{s=1}^{\infty}\dfrac{1}{s(s+1)(s+2)}v^s=\begin{cases}
\dfrac{3v-2}{4v}-\dfrac{(1-v)^2}{2v^2}\log (1-v) & \text{($0<|v|<1$)},\\
                                               0 & \text{($v=0$),}
  \end{cases}\\
L_2(v)=&~\sum_{s=0}^{\infty}\dfrac{1}{s+3}v^{s+1}=\begin{cases}
-\dfrac{1}{v^2}\log (1-v)-\dfrac{2+v}{2v} & \text{($0<|v|<1$)},\\
                                               0 & \text{($v=0$),}
  \end{cases}\\
L_3(v)=&~\sum_{s=0}^{\infty}\dfrac{1}{(s+2)(s+3)}v^{s+1}=\begin{cases}
\dfrac{1-v}{v^2}\log (1-v)+\dfrac{2-v}{2v} & \text{($0<|v|<1$)},\\
                                               0 & \text{($v=0$),}
  \end{cases}
\end{align*}
respectively, and let $B(v)=\sum_{s=0}^{\infty}b_sv^s$. Then, $B(v)$ can be rewritten as 
\begin{align*}
B(v)=&~\dfrac{2}{v\kappa_{\lambda}^{(2)}}\Biggl[L_1\left(\dfrac{n_1-p-\tfrac{1}{2}}{n-p_1-\tfrac{1}{2}}v\right)-L_1\left(\dfrac{n_1-p-\tfrac{1}{2}}{n-\tfrac{1}{2}}v\right)\\
&+\tau_1^2\left\{L_1\left(\tau_1v\right)-L_1\left(\tau_1\dfrac{n_1-p-\tfrac{1}{2}}{n_1-p_1-\tfrac{1}{2}}v\right)\right\}\\
&-\dfrac{1}{n^2}L_2\left(\dfrac{n_1-p-\tfrac{1}{2}}{n}v\right)
-\dfrac{p_1+\tau_1p_2}{n}L_3\left(\dfrac{n_1-p-\tfrac{1}{2}}{n}v\right)\Biggr].
\end{align*}
The difference between $\varphi_{T}(t)$ and $\varphi_{T,s}(t)$ is
\begin{align*}
\varphi_{T}(t)-\varphi_{T,s}(t)
=&~\exp\Bigl(-\dfrac{t^2}{2}\Bigr)\Biggl\{\sum_{k=1}^{s}\dfrac{(it)^{3k}}{k!}\sum_{j=s-k+1}^{\infty}\gamma_{k,j}(it)^j\\
&+\sum_{k=s+1}^{\infty}\dfrac{(it)^{3k}}{k!}\Bigl(\sum_{j=0}^{\infty}\dfrac{\widetilde{\kappa}_{T}^{(j+3)}}{(j+3)!}(it)^j\Bigr)^k\Biggr\}.
\end{align*}
Using \eqref{eqn:a4}, for $|t|\leq mv$, we have
\begin{align*}
\dfrac{1}{|t|}|\varphi_{T}(t)-\varphi_{T,s}(t)|<&~\dfrac{1}{m^{s+1}}\exp\left(-\dfrac{t^2}{2}\right)\Biggl\{\sum_{k=1}^{s}\dfrac{1}{k!}|t|^{s+2k}R_{k,s-k+1}(v)\\
&+\dfrac{|t|^{3s+2}}{(s+1)!}(B(v))^{s+1}\exp\left(t^2vB(v)\right)\Biggr\},
\end{align*}
where
\begin{align*}
R_{k,\ell}(v)=v^{-\ell}\Biggl\{(B(v))^k-\sum_{j=0}^{\ell-1}\Bigl(\sum_{s_1+\cdots+s_k=j}b_{s_1}\cdots b_{s_k}\Bigr)v^j\Biggr\}.
\end{align*}
Hence, we obtain the upper bound $U_1(v)$ for $I_1(v)$:
\begin{align}
\begin{split}\label{eqn:a7}
U_1(v)=&~\dfrac{2}{m^{s+1}}\Biggl\{\sum_{k=1}^{s}\dfrac{1}{k!}
R_{k,s-k+1}(v)\int_{0}^{mv}t^{s+2k}
\exp\left(-\dfrac{t^2}{2}\right)dt\\
&+\dfrac{1}{(s+1)!}(B(v))^{s+1}
\int_{0}^{mv}t^{3s+2}
\exp\left(-\dfrac{t^2}{2}c_v\right)dt\Biggr\},
\end{split}
\end{align}
where $c_v=1-2vB(v)$. If $c_v>0$, $U_1(v)=O(m^{-(s+1)})~(m\rightarrow\infty)$.

We next calculate $I_2(v)$ and derive an upper bound $U_2(v)$ for $I_2(v)$. From \eqref{eqn:a5},
\begin{align}
I_2(v)=2\Biggl\{\int_{mv}^{\infty}\exp\Bigl(-\dfrac{t^2}{2}\Bigr)t^{-1}dt+\sum_{k=1}^{s}\dfrac{1}{k!}\sum_{j=0}^{s-k}\gamma_{k,j}\int_{mv}^{\infty}\exp\Bigl(-\dfrac{t^2}{2}\Bigr)t^{3k+j-1}dt\Biggr\}.\label{eqn:a8}
\end{align}
For $c\in(0,1)$ and $k>-1$, the following inequalities hold.
\begin{align*}
2\int_{mv}^{\infty}\exp\Bigl(-\dfrac{t^2}{2}\Bigr)t^{k}dt&=2\int_{mv}^{\infty}\exp\Bigl(-\dfrac{t^2}{2}(1-c+c)\Bigr)t^{k}dt\\
&<\exp\Bigl(-\dfrac{m^2v^2}{2}(1-c)\Bigr)\Bigl(\dfrac{c}{2}\Bigr)^{-\tfrac{k+1}{2}}\Gamma\Bigl(\dfrac{k+1}{2}\Bigr),\\
2\int_{mv}^{\infty}\exp\Bigl(-\dfrac{t^2}{2}\Bigr)t^{-1}dt&\leq2\exp\Bigl(-\dfrac{m^2v^2}{2}(1-c)\Bigr)\int_{mv\sqrt{\tfrac{c}{2}}}^{\infty}u^{-1}\exp(-u^2)du\\
&\leq2\exp\Bigl(-\dfrac{m^2v^2}{2}(1-c)\Bigr)\dfrac{1}{mv}\sqrt{\dfrac{2}{c}}\int_{mv\sqrt{\tfrac{c}{2}}}^{\infty}\exp(-u^2)du\\
&<2\exp\Bigl(-\dfrac{m^2v^2}{2}(1-c)\Bigr)\dfrac{1}{mv}\sqrt{\dfrac{2}{c}}\int_{0}^{\infty}\exp(-u^2)du\\
&=\exp\Bigl(-\dfrac{m^2v^2}{2}(1-c)\Bigr)\dfrac{1}{mv}\sqrt{\dfrac{2\pi}{c}}.
\end{align*}
Hence, we obtain the upper bound $U_2(v)$ for $I_2(v)$:
\begin{align}
\begin{split}\label{eqn:a9}
U_2(v)=&~\exp\Bigl(-\dfrac{m^2v^2}{2}(1-c)\Bigr)\Biggl\{\dfrac{1}{mv}\sqrt{\dfrac{2\pi}{c}}\\
&+\sum_{k=1}^{s}\dfrac{1}{k!}\sum_{j=0}^{s-k}\gamma_{k,j}\Bigl(\dfrac{c}{2}\Bigr)^{-\tfrac{3k+j}{2}}\Gamma\Bigl(\dfrac{3k+j}{2}\Bigr)\Biggr\}.
\end{split}
\end{align}
We note that 
\begin{align*}
I_2(v)&=O(\exp\left(-m^2v^2(1-c)/2\right))\qquad(m\rightarrow\infty),\\ 
U_2(v)&=O(\exp\left(-m^2v^2(1-c)/2\right))\qquad(m\rightarrow\infty),
\end{align*}
for fixed $v\in(0,1)$ and $c\in(0,1)$.

Moreover, we derive a simpler upper bound $\widetilde{U}_2(v)$ for $I_2(v)$. For sufficiently large $m$, i.e., when $mv>\sqrt{2\pi/c}$, we obtain the simple upper bound $\widetilde{U}_2(v)$ for $I_2(v)$:
\begin{align}
\widetilde{U}_2(v)=\exp\Bigl(-\dfrac{m^2v^2}{2}(1-c)\Bigr)\Biggl\{1+\sum_{k=1}^{s}\dfrac{1}{k!}\sum_{j=0}^{s-k}\gamma_{k,j}\Bigl(\dfrac{c}{2}\Bigr)^{-\tfrac{3k+j}{2}}\Gamma\Bigl(\dfrac{3k+j}{2}\Bigr)\Biggr\}.\label{eqn:a10}
\end{align}
\cite{akita2010high} derived $\widetilde{U}_2(v)$ using the following method. Since the first definite integral in \eqref{eqn:a8} diverges as $mv\rightarrow+0$, a sufficient condition is required. If $mv\ge\sqrt{2}$,
\begin{align*}
2\int_{mv}^{\infty}\exp\Bigl(-\dfrac{t^2}{2}\Bigr)t^{-1}dt&<\int_{mv}^{\infty}\exp\Bigl(-\dfrac{t^2}{2}\Bigr)tdt\\
&=\exp\Bigl(-\dfrac{m^2v^2}{2}\Bigr)<\exp\Bigl(-\dfrac{m^2v^2}{2}(1-c)\Bigr).
\end{align*}
Therefore, \cite{akita2010high} obtained \eqref{eqn:a10}. For fixed $v\in(0,1)$ and $c\in(0,1)$, we note that
$\widetilde{U}_2(v)=O(\exp\left(-m^2v^2(1-c)/2\right))$~($m\rightarrow\infty$).

Finally, we derive an upper bound $U_3(v)$ for $I_3(v)$. The characteristic function of $T$ can be expressed as
\begin{align*}
\varphi_{T}(t)&=\mathrm{E}[\exp(itT)]\\
&=\exp(-i\tilde t\kappa_{\lambda}^{(1)})\mathrm{E}\Biggl[\exp\Bigl\{i\tilde t\bigl(-\tfrac{2}{N}\log\lambda\bigr)\Bigr\}\Biggr]=\exp(-i\tilde t\kappa_{\lambda}^{(1)})\varphi_{\lambda}(\mkern 1.5mu\tilde t\mkern 1.5mu),
\end{align*}
where $\tilde t=t/(\kappa_{\lambda}^{(2)})^{\tfrac{1}{2}}$. From $\tilde t\kappa_{\lambda}^{(1)}\in\mathbb{R}$, we obtain the following expression by applying Euler's formula.
\begin{align*}
\left|\exp(-i\tilde t\kappa_{\lambda}^{(1)})\right|=\left|\cos(\tilde t\kappa_{\lambda}^{(1)})+i\sin(\tilde t\kappa_{\lambda}^{(1)})\right|=1.\\
\end{align*}
Therefore,
\begin{align*}
\left|\varphi_{T}(t)\right|
=&\left|\varphi_{\lambda}(\mkern 1.5mu\tilde t\mkern 1.5mu)\right|\\
=&\left|K^{-i\tilde t}\right|\Biggl|\prod_{\ell=1}^{p_1}\dfrac{\Gamma\left[\dfrac{1}{2}(n-p_1+\ell)-i\tilde t\right]}{\Gamma\left[\dfrac{1}{2}(n-p_1+\ell)\right]}\prod_{\ell=1}^{p_2}\dfrac{\Gamma\left[\dfrac{1}{2}(n_1-p+\ell)-\tau_1i\tilde t\right]}{\Gamma\left[\dfrac{1}{2}(n_1-p+\ell)\right]}\\
&\times\dfrac{\Gamma\left[\dfrac{1}{2}(np_1+n_1p_2)\right]}{\Gamma\left[\dfrac{1}{2}(np_1+n_1p_2)-(p_1+\tau_1p_2)i\tilde t\right]}\Biggr|,
\end{align*}
where
\begin{align*}
K=\dfrac{(p_1+\tau_1p_2)^{(p_1+\tau_1p_2)}}{\tau_1^{\tau_1p_2}}>0.
\end{align*}
From $\tilde t\log K\in\mathbb{R}$,
\begin{align*}
|K^{-i\tilde t}|&=\left|\exp(-i\tilde t\log K)\right|=\left|\cos(\tilde t\log K)+i\sin(\tilde t\log K)\right|=1.
\end{align*}
Hence, we have
\begin{align*}
\left|\varphi_{T}(t)\right|=&~\Biggl|\prod_{\ell=1}^{p_1}\dfrac{\Gamma\left[\dfrac{1}{2}(n-p_1+\ell)-i\tilde t\right]}{\Gamma\left[\dfrac{1}{2}(n-p_1+\ell)\right]}\prod_{\ell=1}^{p_2}\dfrac{\Gamma\left[\dfrac{1}{2}(n_1-p+\ell)-\tau_1i\tilde t\right]}{\Gamma\left[\dfrac{1}{2}(n_1-p+\ell)\right]}\\
&\times\dfrac{\Gamma\left[\dfrac{1}{2}(np_1+n_1p_2)\right]}{\Gamma\left[\dfrac{1}{2}(np_1+n_1p_2)-(p_1+\tau_1p_2)i\tilde t\right]}\Biggr|.
\end{align*}
It is known that
\begin{align*}
\left|\dfrac{\Gamma[x+iy]}{\Gamma[x]}\right|^2=\prod_{k=0}^{\infty}\Biggl\{1+\dfrac{y^2}{(x+k)^2}\Biggr\}^{-1},
\end{align*}
for any $x,y\in \mathbb{R}~(x>0)$. Then,
\begin{align*}
\log|\varphi_{T}(t)|=&~\dfrac{1}{2}\sum_{k=0}^{\infty}\log\left\{1+\dfrac{(p_1+\tau_1p_2)^2{\tilde t}^{\,2}}{\Bigl(\tfrac{1}{2}(np_1+n_1p_2)+k\Bigr)^2}\right\}\\
&-\dfrac{1}{2}\sum_{\ell=1}^{p_1}\sum_{k=0}^{\infty}\log\left\{1+\dfrac{{\tilde t}^{\,2}}{\Bigl(\tfrac{1}{2}(n-p_1+\ell)+k\Bigr)^2}\right\}\\
&-\dfrac{1}{2}\sum_{\ell=1}^{p_2}\sum_{k=0}^{\infty}\log\left\{1+\dfrac{\tau_1^2{\tilde t}^{\,2}}{\Bigl(\tfrac{1}{2}(n_1-p+\ell)+k\Bigr)^2}\right\}.
\end{align*}
Using 
\begin{align*}
\log\left(1+\dfrac{a}{x^2}\right)=\displaystyle\int_{0}^{a}\dfrac{1}{x^2+s}ds,
\end{align*}
for any $a>0$ and $x>0$,
\begin{align}
\begin{split}\label{eqn:a11}
\log|\varphi_{T}(t)|=&~\dfrac{1}{2}\sum_{k=0}^{\infty}\displaystyle\int_{0}^{(p_1+\tau_1p_2)^2{\tilde t}^{\,2}}\dfrac{1}{\Bigl(\tfrac{1}{2}(np_1+n_1p_2)+k\Bigr)^2+s}ds\\
&-\dfrac{1}{2}\sum_{\ell=1}^{p_1}\sum_{k=0}^{\infty}\displaystyle\int_{0}^{{\tilde t}^{\,2}}\dfrac{1}{\Bigl(\tfrac{1}{2}(n-p_1+\ell)+k\Bigr)^2+s}ds\\
&-\dfrac{1}{2}\sum_{\ell=1}^{p_2}\sum_{k=0}^{\infty}\displaystyle\int_{0}^{\tau_1^2{\tilde t}^{\,2}}\dfrac{1}{\Bigl(\tfrac{1}{2}(n_1-p+\ell)+k\Bigr)^2+s}ds.
\end{split}
\end{align}
Here, for any $a>0$ and $s>0$, setting $f_k(s)=1/\{(a+k)^2+s\}$ ensures the uniform convergence of $f_k(s)$. Therefore, in \eqref{eqn:a11}, the infinite series can be interchanged with the definite integral as
\begin{align}
\begin{split}\label{eqn:a12}
\log|\varphi_{T}(t)|=&~\dfrac{1}{2}\displaystyle\int_{0}^{{\tilde t}^{\,2}}\Biggl\{(p_1+\tau_1p_2)^2\sum_{k=0}^{\infty}\dfrac{1}{\Bigl(\tfrac{1}{2}(np_1+n_1p_2)+k\Bigr)^2+(p_1+\tau_1p_2)^2u}\\
&-\sum_{\ell=1}^{p_1}\sum_{k=0}^{\infty}\dfrac{1}{\Bigl(\tfrac{1}{2}(n-p_1+\ell)+k\Bigr)^2+u}\\
&-\tau_1^2\sum_{\ell=1}^{p_2}\sum_{k=0}^{\infty}\dfrac{1}{\Bigl(\tfrac{1}{2}(n_1-p+\ell)+k\Bigr)^2+\tau_1^2u}\Biggr\}du.
\end{split}
\end{align}
Here, let $g(x)>0$ be a monotonically decreasing function on $x>0$. For $a>0$, we use the following inequality. 
\begin{align}
\displaystyle\int_{0}^{\infty}g(a+x)dx<\sum_{k=0}^{\infty}g(a+k)<g(a)+\displaystyle\int_{0}^{\infty}g(a+x)dx.\label{eqn:a13}
\end{align}
In \eqref{eqn:a12}, $g(x)$ is the monotonically decreasing function and we obtain the following inequality for $\omega>0$.
\begin{align}
\begin{split}\label{eqn:a14}
\dfrac{1}{\sqrt{\omega}}\left(\dfrac{\pi}{2}-\arctan\dfrac{a}{\sqrt{\omega}}\right)&<\sum_{k=0}^{\infty}\dfrac{1}{(a+k)^2+\omega}\\
&<\dfrac{1}{a^2+\omega}+\dfrac{1}{\sqrt{\omega}}\left(\dfrac{\pi}{2}-\arctan\dfrac{a}{\sqrt{\omega}}\right).
\end{split}
\end{align}
Using the inequality \eqref{eqn:a14}, \eqref{eqn:a12} can be bounded as
\begin{align*}
\log|\varphi_{T}(t)|<&~\dfrac{1}{2}\displaystyle\int_{0}^{{\tilde t}^{\,2}}\dfrac{1}{\sqrt{u}}\Biggl[\sum_{\ell=1}^{p_1}\arctan{\left\{\dfrac{1}{\sqrt{u}}\left(\dfrac{1}{2}(n-p_1+\ell)\right)\right\}}\\
&+\tau_1\sum_{\ell=1}^{p_2}\arctan{\left\{\dfrac{1}{\tau_1\sqrt{u}}\left(\dfrac{1}{2}(n_1-p+\ell)\right)\right\}}\\
&-(p_1+\tau_1p_2)\arctan{\left\{\dfrac{1}{(p_1+\tau_1p_2)\sqrt{u}}\left(\dfrac{1}{2}(np_1+n_1p_2)\right)\right\}}\Biggr]\\
&+\dfrac{(p_1+\tau_1p_2)^2}{\left\{\dfrac{1}{2}(np_1+n_1p_2)\right\}^2+(p_1+\tau_1p_2)^2u}du\\
=&~\displaystyle\int_{0}^{{\tilde t}}\Biggl[\sum_{\ell=1}^{p_1}\arctan{\left\{\dfrac{1}{y}\left(\dfrac{1}{2}(n-p_1+\ell)\right)\right\}}\\
&+\tau_1\sum_{\ell=1}^{p_2}\arctan{\left\{\dfrac{1}{\tau_1y}\left(\dfrac{1}{2}(n_1-p+\ell)\right)\right\}}\\
&-(p_1+\tau_1p_2)\arctan{\left\{\dfrac{1}{(p_1+\tau_1p_2)y}\left(\dfrac{1}{2}(np_1+n_1p_2)\right)\right\}}\Biggr]dy\\
&+\dfrac{1}{2}\log\left\{\dfrac{1}{4}(np_1+n_1p_2)^2+(p_1+\tau_1p_2)^2{\tilde t}^{\,2}\right\}-\log\left\{\dfrac{1}{2}(np_1+n_1p_2)\right\}\\[3pt]
=&-\tilde{t}^{\,2}F(\tilde t; n,n_1,p_1,p_2),
\end{align*}
where
\begin{align*}
F(t;n,n_1,p_1,p_2)
=&~\dfrac{1}{t^2}\displaystyle\int_{0}^{t}\Biggl[(p_1+\tau_1p_2)\arctan{\left\{\dfrac{1}{(p_1+\tau_1p_2)y}\left(\dfrac{1}{2}(np_1+n_1p_2)\right)\right\}}\\
&-\sum_{\ell=1}^{p_1}\arctan{\left\{\dfrac{1}{y}\left(\dfrac{1}{2}(n-p_1+\ell)\right)\right\}}\\
&-\tau_1\sum_{\ell=1}^{p_2}\arctan{\left\{\dfrac{1}{\tau_1y}\left(\dfrac{1}{2}(n_1-p+\ell)\right)\right\}}\Biggr]dy\\
&+\dfrac{1}{t^2}\log\left\{\dfrac{1}{2}(np_1+n_1p_2)\right\}\\
&-\dfrac{1}{2t^2}\log\left\{\dfrac{1}{4}(np_1+n_1p_2)^2+(p_1+\tau_1p_2)^2t^2\right\}.
\end{align*}
By applying a change of variables, $F(t;n,n_1,p_1,p_2)$ can be expressed as
\begin{align}
\begin{split}\label{eqn:a15}
F(t;n,n_1,p_1,p_2)=&~\dfrac{1}{t^2}\Biggl\{\dfrac{1}{2}(np_1+n_1p_2)H\left[\dfrac{np_1+n_1p_2}{2(p_1+\tau_1p_2)t}\right]\\
&+\log\left(\dfrac{1}{2}(np_1+n_1p_2)\right)\\
&-\dfrac{1}{2}\sum_{\ell=1}^{p_1}(n-p_1+\ell)H\left[\dfrac{n-p_1+\ell}{2t}\right]\\
&-\dfrac{1}{2}\sum_{\ell=1}^{p_2}(n_1-p+\ell)H\left[\dfrac{n_1-p+\ell}{2\tau_1t}\right]\\
&-\dfrac{1}{2}\log\left(\dfrac{1}{4}(np_1+n_1p_2)^2+(p_1+\tau_1p_2)^2t^2\right)\Biggr\},
\end{split}
\end{align}
where
\begin{align*}
H(z)=\int_{z}^{\infty}\dfrac{1}{u^2}\arctan{u}du=\dfrac{1}{z}\arctan{z}-\log z+\dfrac{1}{2}\log(1+z^2),
\end{align*}
for $z>0$. Therefore,
\begin{align*}
|\varphi_{T}(t)|<\exp\left\{-{\tilde t}^{\,2}F\left({\tilde t};n,n_1,p_1,p_2\right)\right\}.
\end{align*}
Hence, we obtain the upper bound $U_3(v)$ for $I_3(v)$:
\begin{align}
U_3(v)=2\int_{m_0v}^{\infty}\dfrac{1}{t}\exp\left\{-t^2F\left(t;n,n_1,p_1,p_2\right)\right\}dt,\label{eqn:a16}
\end{align}
where
\begin{align*}
m_0=\dfrac{n_1-p-\tfrac{1}{2}}{2}.
\end{align*}
$U_3(v)$ diverges if $p<3$, whereas $U_3(v)=O(m^{-(p^2-p-4)/4})~(m\rightarrow\infty)$ if $p\ge3$.
The results obtained here are summarized in the following theorem.
\begin{thm}\label{thm:2}
Let $Q_s(x)$ be the Edgeworth expansion for the null distribution of the LRT statistic up to the order $O(m^{-s})$ given by \eqref{eqn:a6} under high-dimensional two-step monotone incomplete data. Then, 
\begin{align*}
\sup_{x}\bigl|\Pr(T\leq x)-Q_s(x)\bigr|
<\dfrac{1}{2\pi}\Bigl(U_1(v)+I_2(v)+U_3(v)\Bigr),
\end{align*}
for $v\in(0,1)$ and $c\in(0,1)$, where $U_1(v)$ is given by \eqref{eqn:a7}, $I_2(v)$ is given by \eqref{eqn:a8}, and $U_3(v)$ is given by \eqref{eqn:a16} with $F\left(t;n,n_1,p_1,p_2\right)$ given by \eqref{eqn:a15}. Two upper bounds for $I_2(v)$ are also obtained: $U_2(v)$ and a simpler one, $\widetilde{U}_2(v)$, where $U_2(v)$ is given by \eqref{eqn:a9}, and $\widetilde{U}_2(v)$ is given by \eqref{eqn:a10}.
\end{thm}
\section{Numerical experiments}
\label{sec:Numerical experiments}
This section presents numerical results to evaluate the Edgeworth expansion $Q_2(x)$ and its computable error bounds. Specifically, Section \ref{subsec:Type I error rates} employs Monte Carlo simulations to compare the approximation accuracy of the Edgeworth expansion $Q_2(x)$ and the large-sample asymptotic expansion for the null distribution of the LRT statistic. In Section \ref{subsec:Example}, we evaluate the MAE between the distribution function of the standardized test statistic $T$ and the Edgeworth expansion $Q_2(x)$ via Monte Carlo simulations. Furthermore, numerical experiments are conducted to analyze the behavior of four computable error bounds. Finally, in Section \ref{subsec:Cornish-Fisher}, we derive a high-dimensional Cornish-Fisher expansion to evaluate whether the empirical Type I error rate is precisely controlled, and we also assess the empirical power.
\subsection{Comparison of approximation accuracy}
\label{subsec:Type I error rates}
\begin{figure}[H]
\centering
\includegraphics[width=0.9\columnwidth]{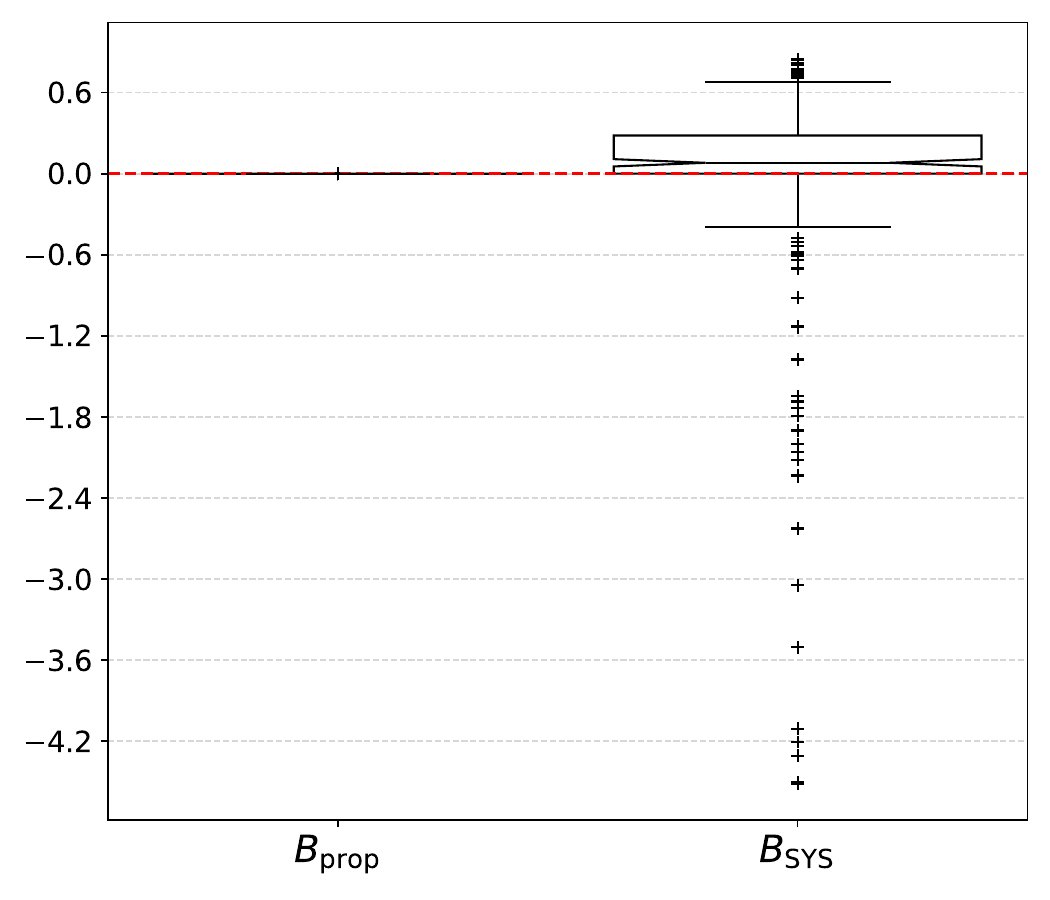}
\caption{Box plot of the biases $B_{\mathrm{prop}}$ and $B_{\mathrm{SYS}}$. The red dashed line indicates the case where the bias is zero.}
\label{fig:1}
\end{figure}

\begin{figure}[H]
\centering
\includegraphics[width=0.89\columnwidth]{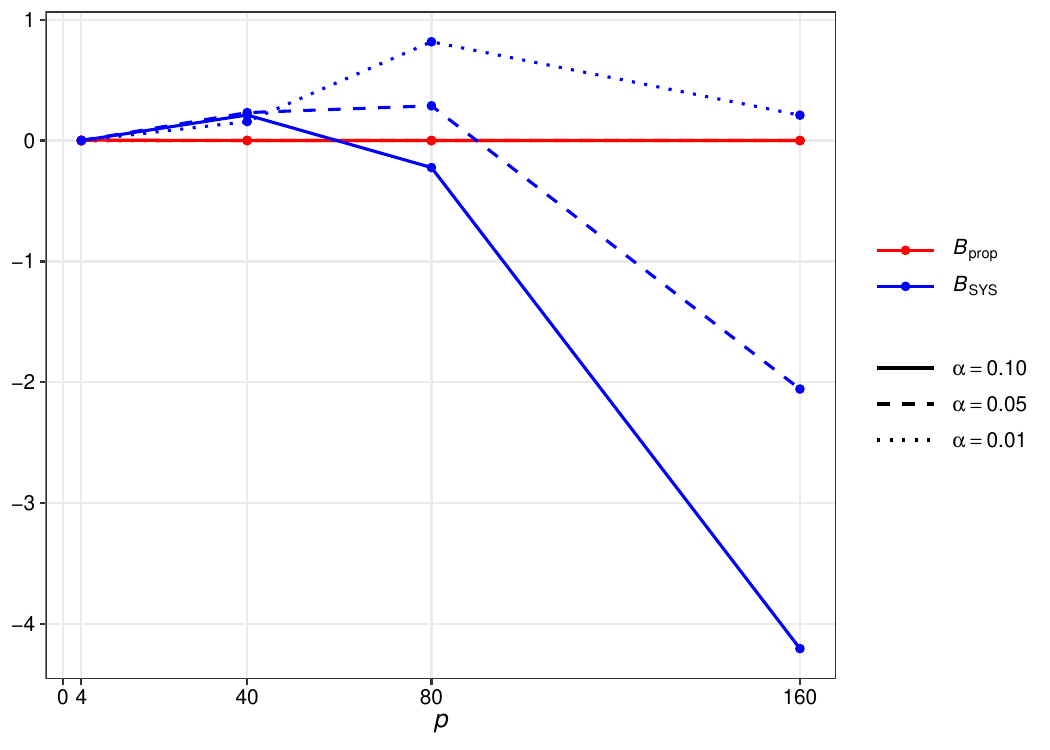}
\caption{Line graph of the biases $B_{\mathrm{prop}}$ and $B_{\mathrm{SYS}}$ for $N_1=N_2=200$ and $p_1=p_2=2,20,40,80$. The red lines represent $B_{\mathrm{prop}}$ and blue lines represent $B_{\mathrm{SYS}}$. The solid, dashed, and dotted lines indicate $\alpha=0.10,0.05,$ and $0.01$, respectively.}
\label{fig:2}
\end{figure}

We compare the approximation accuracy of the Edgeworth expansion $Q_2(x)$ for $s=2$ and the large-sample asymptotic expansion for the null distribution of the LRT statistic, given by \eqref{eqn:a6} and \eqref{eqn:a2}, respectively, using $10^6$ Monte Carlo simulations. We assess which approximation provides a closer match to the Type I error rate in both large-sample and high-dimensional settings. Let the Type I error rate under two-step monotone incomplete data be as follows.
\begin{align*}
\alpha_1&=\mathrm{Pr}\{-2\log\lambda>\chi_f^{2}(\alpha)\},
\end{align*}
where $f=(p+2)(p-1)/2$ and $\chi_k^{2}(\alpha)$ is the upper $100\alpha$ percentile of the $\chi^2$ distribution with $k$ degrees of freedom. From \eqref{eqn:a6}, $Q_2(x)$ can be written as
\begin{align*}
Q_2(x)
=&~\Phi(x)-\phi(x)\Biggl\{\dfrac{1}{6}\widetilde{\kappa}_{T}^{(3)}(x^2-1)+\dfrac{1}{24}\widetilde{\kappa}_{T}^{(4)}(x^3-3x)\\
&+\dfrac{1}{72}(\widetilde{\kappa}_{T}^{(3)})^2(x^5-10x^3+15x)\Biggr\}.
\end{align*}
$\alpha_1$ can also be rewritten as $\alpha_1=\mathrm{Pr}\{T>q(\alpha)\},$
where 
\begin{align*}
q(\alpha)=\dfrac{\tfrac{1}{N}\chi_f^{2}(\alpha)-\kappa_{\lambda}^{(1)}}{(\kappa_{\lambda}^{(2)})^{\tfrac{1}{2}}}.
\end{align*}
In this case, we compare the biases $B_{\mathrm{prop}}$ and $B_{\mathrm{SYS}}$, where
\begin{align*}
B_{\mathrm{prop}}=&~\alpha_1-\{1-Q_2(q(\alpha))\},\\
B_{\mathrm{SYS}}=&~\alpha_1-\Biggl[1-\Bigl[G_f(\chi_f^{2}(\alpha))+\dfrac{\beta}{N}\bigl\{G_{f+2}(\chi_f^{2}(\alpha))-G_{f}(\chi_f^{2}(\alpha))\bigr\}\\
&+\dfrac{\gamma}{N^2}\bigl\{G_{f+4}(\chi_f^{2}(\alpha))-G_{f}(\chi_f^{2}(\alpha))\bigr\}\Bigr]\Biggr].
\end{align*}
$\beta$ and $\gamma$ are given by Theorem \ref{thm:1}.
In the simulations, we generate samples with $\bm{\mu}=\bm{0}$ and $\bm{\Sigma}=\bm{I}_p$. The configuration of the values of the sample sizes and dimensions are $(i)~N_1=50,100,200$, $N_1/N_2=0.5,1,2$, and $p_1=p_2=2$; and $(ii)~N_1=50,100,200$, $N_1/N_2=0.5,1,2$, $p/N_1=0.2,0.4,0.8$, and $p_1/p_2=0.25,1,4$. Furthermore, the nominal significance levels are $\alpha=0.10,0.05,0.01$. Figure \ref{fig:1} shows the box plot of the biases $B_{\mathrm{prop}}$ and $B_{\mathrm{SYS}}$ for all configurations of $N_1,N_2,p_1,p_2$, and $\alpha$.
As shown in Figure \ref{fig:1}, $B_{\mathrm{prop}}$ outperforms better than $B_{\mathrm{SYS}}$, tending to yield smaller bias across configurations. 
Figure \ref{fig:2} presents the line graph of the biases $B_{\mathrm{prop}}$ and $B_{\mathrm{SYS}}$ for $N_1=N_2=200$ and $p_1=p_2=2,20,40,80$. Figure \ref{fig:2} indicates that $B_{\mathrm{prop}}$ remains close to zero, whereas $B_{\mathrm{SYS}}$ tends to take values increasingly farther from zero with growing $p$. The results for each configuration of $N_1,N_2,p_1,p_2$, and $\alpha$ are given in \ref{app:D}.
\subsection{Evaluation of error bounds and maximum absolute error}
\label{subsec:Example}
We evaluate the error bounds and the MAE between the distribution function of $T$ and the Edgeworth expansion $Q_2(x)$ for $s=2$, using Monte Carlo simulation with $10^6$ replications. The settings for $\bm\mu$ and $\bm\Sigma$ are identical to those in Section \ref{subsec:Type I error rates}. $U_1(v)$ can be expressed as
\begin{align*}
U_1(v)=&~\dfrac{2}{m^{3}}\int_{0}^{mv}\Biggl\{\left(R_{1,2}(v)t^4+\dfrac{1}{2}R_{2,1}(v)t^6\right)\exp\left(-\dfrac{t^2}{2}\right)\\
&+\dfrac{1}{6}(B(v))^{3}t^{8}\exp\left(-\dfrac{t^2}{2}c_v\right)\Biggr\}dt,
\end{align*}
where
\begin{align*}
R_{1,2}(v)=\dfrac{1}{v^2}(B(v)-b_0-b_1v),~R_{2,1}(v)=\dfrac{1}{v}\left\{(B(v))^2-b_0^2\right\}.
\end{align*}
$I_2(v)$, $U_2(v)$ and $\widetilde{U}_2(v)$ can be written as
\begin{align*}
I_2(v)=&~2\int_{mv}^{\infty}\left\{\dfrac{1}{t}+\dfrac{1}{6}\widetilde{\kappa}_{T}^{(3)}t^2+\dfrac{1}{72}(\widetilde{\kappa}_{T}^{(3)})^2t^5+\dfrac{1}{24}\widetilde{\kappa}_{T}^{(4)}t^3\right\}\exp\left(-\dfrac{t^2}{2}\right)dt,\\
U_2(v)=&~\exp\Bigl(-\dfrac{m^2v^2}{2}(1-c)\Bigr)\Biggl\{\dfrac{1}{mv}\sqrt{\dfrac{2\pi}{c}}\\
&+\gamma_{1,0}\left(\dfrac{c}{2}\right)^{-\tfrac{3}{2}}\Gamma\left(\dfrac{3}{2}\right)+\gamma_{1,1}\Bigl(\dfrac{c}{2}\Bigr)^{-2}\Gamma(2)+\gamma_{2,0}\Bigl(\dfrac{c}{2}\Bigr)^{-3}\Gamma(3)\Biggr\},\\
\widetilde{U}_2(v)=&~\exp\Bigl(-\dfrac{m^2v^2}{2}(1-c)\Bigr)\Biggl\{1+\gamma_{1,0}\left(\dfrac{c}{2}\right)^{-\tfrac{3}{2}}\Gamma\left(\dfrac{3}{2}\right)\\
&+\gamma_{1,1}\Bigl(\dfrac{c}{2}\Bigr)^{-2}\Gamma(2)+\gamma_{2,0}\Bigl(\dfrac{c}{2}\Bigr)^{-3}\Gamma(3)\Biggr\},
\end{align*}
respectively, where $\gamma_{1,0}=\widetilde{\kappa}_{T}^{(3)}/6$, $\gamma_{1,1}=\widetilde{\kappa}_{T}^{(4)}/24$, $\gamma_{2,0}=(\widetilde{\kappa}_{T}^{(3)})^2/36$.
$U_3(v)$ does not depend on $s$. Although the upper bounds in Theorem \ref{thm:2} can be minimized numerically with respect to $v$ and $c$, in practice it suffices to evaluate the error bounds on the grids $v=0.05,0.10,\dots,0.95$, and $c=0.05,0.10,\dots,0.95$, and then choose the minimum. The MAE and the four kinds of error bounds are given by $\mathrm{MAE}=\sup_{x}\bigl|\Pr(T\leq x)-Q_2(x)\bigr|$ and
\begin{align*}
\mathrm{BOUND1}&=\min_{\substack{v=0.05,0.10,\dots,0.95}}\dfrac{1}{2\pi}\Bigl(U_1(v)+I_2(v)+U_3(v)\Bigr),\\
\mathrm{BOUND2}&=\min_{\substack{v=0.05,0.10,\dots,0.95\\
c=0.05,0.10,\dots,0.95}}\dfrac{1}{2\pi}\Bigl(U_1(v)+U_2(v)+U_3(v)\Bigr),\\
\mathrm{BOUND3}&=\min_{\substack{v=0.05,0.10,\dots,0.95\\
c=0.05,0.10,\dots,0.95,~c>\tfrac{2\pi}{m^2v^2}}}\dfrac{1}{2\pi}\left(U_1(v)+\widetilde{U}_2(v)+U_3(v)\right),\\
\mathrm{BOUND4}&=\min_{\substack{v=0.05,0.10,\dots,0.95,~v\ge\tfrac{\sqrt2}{m}\\
c=0.05,0.10,\dots,0.95}}\dfrac{1}{2\pi}\left(U_1(v)+\widetilde{U}_2(v)+U_3(v)\right).
\end{align*}
$\mathrm{BOUND1}$ is the upper bound obtained from the improper integral $I_2(v)$. $\mathrm{BOUND2}$ is the upper bound derived from $U_2(v)$, the upper bound on $I_2(v)$ obtained in this work. $\mathrm{BOUND3}$ and $\mathrm{BOUND4}$ are based on $\widetilde{U}_2(v)$, the simpler upper bound for $I_2(v)$ given in \cite{akita2010high}. However, we note that $\mathrm{BOUND3}$ and $\mathrm{BOUND4}$ rely on different sufficient conditions. Let $v_{\min}$ and $c_{\min}$ denote the values of $v$ and $c$ that minimize each of the error bounds $\mathrm{BOUND1}$--$\mathrm{BOUND4}$, and let $c_{v_{\min}}$ denote $c_v$ evaluated at $v=v_{\min}$. If $c_{v_{\min}}>0$, then the error bounds converge to zero as $m\rightarrow\infty$.

We begin by evaluating the inequality stated in Theorem \ref{thm:2}. Table \ref{tab:1} shows that the MAE is smaller than each of $\mathrm{BOUND1}$--$\mathrm{BOUND4}$.

We then evaluate the error bounds $\mathrm{BOUND1}$--$\mathrm{BOUND4}$. Tables \ref{tab:1}--\ref{tab:3} report results with $n$ and $n_1$ fixed; Table \ref{tab:4} reports the results with $p$ and $p_1$ fixed; and Table \ref{tab:5} reports the results when $n,n_1,p$, and $p_1$ increase simultaneously. From Tables \ref{tab:1}--\ref{tab:5}, we observe $c_{v_{min}}>0$. When both $p$ and $n_1-p$ are moderately large, the error bounds are small enough for practical use. Moreover, as $n,n_1,p$, and $p_1$ increase, $m$ grows and all error bounds converge to zero.
\begin{table}[H]
\centering
\caption{$\mathrm{MAE}$ and error bounds for $n=60$ and $n_1=40$. The values in parentheses represent the parameters used to minimize each bound. For $\mathrm{BOUND1}$, the parameters are $(v_{\min}, c_{v_{\min}})$. For $\mathrm{BOUND2}$, $\mathrm{BOUND3}$, and $\mathrm{BOUND4}$, the parameters are $(v_{\min}, c_{\min}, c_{v_{\min}})$.}
\label{tab:1}
\footnotesize
\setlength{\tabcolsep}{3pt} 
\renewcommand{\arraystretch}{1.2}
\resizebox{\textwidth}{!}{
\begin{tabular}{rrccccccc}
\toprule
$p_1$ & $p_2$ & $\mathrm{MAE}$ & $\mathrm{BOUND1}$ & $\mathrm{BOUND2}$ & $\mathrm{BOUND3}$ & $\mathrm{BOUND4}$ & $\kappa_{\lambda}^{(2)}$ & $m$ \\
\midrule
5 & 5 & 0.0012 & 0.0287 (0.85, 0.50) & 0.0943 (0.95, 0.40, 0.39) & 0.1796 (0.95, 0.90, 0.39) & 0.0856 (0.95, 0.40, 0.39) & 0.037 & 2.83 \\
5 & 10 & 0.0008 & 0.0087 (0.75, 0.65) & 0.0242 (0.90, 0.30, 0.54) & 0.0348 (0.95, 0.55, 0.50) & 0.0230 (0.90, 0.30, 0.54) & 0.092 & 3.71 \\
5 & 15 & 0.0011 & 0.0036 (0.75, 0.71) & 0.0082 (0.90, 0.25, 0.62) & 0.0100 (0.95, 0.40, 0.59) & 0.0080 (0.90, 0.25, 0.62) & 0.185 & 4.19 \\
\addlinespace[0.5em]
10 & 5 & 0.0006 & 0.0086 (0.80, 0.62) & 0.0243 (0.90, 0.30, 0.55) & 0.0337 (0.95, 0.55, 0.51) & 0.0229 (0.90, 0.30, 0.55) & 0.089 & 3.65 \\
10 & 10 & 0.0008 & 0.0036 (0.75, 0.71) & 0.0083 (0.90, 0.25, 0.62) & 0.0107 (0.90, 0.45, 0.62) & 0.0080 (0.90, 0.25, 0.62) & 0.182 & 4.16 \\
10 & 15 & 0.0009 & 0.0020 (0.80, 0.74) & 0.0040 (0.95, 0.20, 0.67) & 0.0049 (0.95, 0.40, 0.67) & 0.0039 (0.95, 0.20, 0.67) & 0.336 & 4.20 \\
\addlinespace[0.5em]
15 & 5 & 0.0007 & 0.0037 (0.80, 0.69) & 0.0086 (0.90, 0.25, 0.64) & 0.0107 (0.95, 0.45, 0.61) & 0.0082 (0.90, 0.25, 0.64) & 0.173 & 4.06 \\
15 & 10 & 0.0007 & 0.0021 (0.80, 0.75) & 0.0042 (0.95, 0.20, 0.68) & 0.0060 (0.95, 0.45, 0.68) & 0.0040 (0.95, 0.20, 0.68) & 0.328 & 4.15 \\
15 & 15 & 0.0010 & 0.0017 (0.95, 0.75) & 0.0062 (0.95, 0.25, 0.75) & 0.0145 (0.95, 0.55, 0.75) & 0.0054 (0.95, 0.25, 0.75) & 0.596 & 3.67 \\
\bottomrule
\end{tabular}
}
\end{table}
\begin{table}[H]
\centering
\caption{Error bounds for $n=60$ and $n_1=50$. The values in parentheses represent the parameters used to minimize each bound. For $\mathrm{BOUND1}$, the parameters are $(v_{\min}, c_{v_{\min}})$. For $\mathrm{BOUND2}$, $\mathrm{BOUND3}$, and $\mathrm{BOUND4}$, the parameters are $(v_{\min}, c_{\min}, c_{v_{\min}})$.}
\label{tab:2}
\footnotesize
\setlength{\tabcolsep}{4pt} 
\renewcommand{\arraystretch}{1.2}
\resizebox{\textwidth}{!}{%
\begin{tabular}{rrcccccc}
\toprule
$p_1$ & $p_2$ & $\mathrm{BOUND1}$ & $\mathrm{BOUND2}$ & $\mathrm{BOUND3}$ & $\mathrm{BOUND4}$ & $\kappa_{\lambda}^{(2)}$ & $m$ \\
\midrule
5 & 5 & 0.0279 (0.65, 0.51) & 0.0914 (0.75, 0.40, 0.37) & 0.1535 (0.80, 0.75, 0.29) & 0.0847 (0.75, 0.40, 0.37) & 0.035 & 3.70 \\
5 & 10 & 0.0083 (0.60, 0.62) & 0.0233 (0.65, 0.30, 0.57) & 0.0320 (0.70, 0.55, 0.51) & 0.0218 (0.65, 0.30, 0.57) & 0.084 & 5.01 \\
5 & 15 & 0.0032 (0.55, 0.70) & 0.0072 (0.65, 0.25, 0.62) & 0.0092 (0.65, 0.45, 0.62) & 0.0070 (0.65, 0.25, 0.62) & 0.163 & 5.96 \\
\addlinespace[0.5em]
10 & 5 & 0.0082 (0.60, 0.62) & 0.0234 (0.70, 0.30, 0.52) & 0.0315 (0.70, 0.55, 0.52) & 0.0218 (0.65, 0.30, 0.57) & 0.083 & 4.98 \\
10 & 10 & 0.0032 (0.55, 0.70) & 0.0072 (0.65, 0.25, 0.63) & 0.0092 (0.65, 0.45, 0.63) & 0.0070 (0.65, 0.25, 0.63) & 0.162 & 5.94 \\
10 & 15 & 0.0016 (0.55, 0.74) & 0.0030 (0.65, 0.20, 0.68) & 0.0036 (0.65, 0.40, 0.68) & 0.0029 (0.65, 0.20, 0.68) & 0.282 & 6.51 \\
\addlinespace[0.5em]
15 & 5 & 0.0032 (0.55, 0.71) & 0.0071 (0.65, 0.25, 0.63) & 0.0092 (0.70, 0.40, 0.59) & 0.0069 (0.65, 0.25, 0.63) & 0.159 & 5.89 \\
15 & 10 & 0.0016 (0.55, 0.74) & 0.0030 (0.65, 0.20, 0.68) & 0.0036 (0.65, 0.40, 0.68) & 0.0029 (0.65, 0.20, 0.68) & 0.279 & 6.47 \\
15 & 15 & 0.0010 (0.55, 0.78) & 0.0016 (0.65, 0.20, 0.73) & 0.0017 (0.70, 0.30, 0.70) & 0.0016 (0.65, 0.20, 0.73) & 0.458 & 6.60 \\
\bottomrule
\end{tabular}%
}
\end{table}
\begin{table}[H]
\centering
\caption{Error bounds for $n=120$ and $n_1=80$. The values in parentheses represent the parameters used to minimize each bound. For $\mathrm{BOUND1}$, the parameters are $(v_{\min}, c_{v_{\min}})$. For $\mathrm{BOUND2}$, $\mathrm{BOUND3}$, and $\mathrm{BOUND4}$, the parameters are $(v_{\min}, c_{\min}, c_{v_{\min}})$.}
\label{tab:3}
\footnotesize
\setlength{\tabcolsep}{4pt} 
\renewcommand{\arraystretch}{1.2}
\resizebox{\textwidth}{!}{%
\begin{tabular}{rrcccccc}
\toprule
$p_1$ & $p_2$ & $\mathrm{BOUND1}$ & $\mathrm{BOUND2}$ & $\mathrm{BOUND3}$ & $\mathrm{BOUND4}$ & $\kappa_{\lambda}^{(2)}$ & $m$ \\
\midrule
5 & 5 & 0.0282 (0.75, 0.52) & 0.0921 (0.85, 0.40, 0.41) & 0.1489 (0.90, 0.80, 0.34) & 0.0834 (0.85, 0.40, 0.41) & 0.008 & 3.16 \\
5 & 10 & 0.0079 (0.65, 0.64) & 0.0222 (0.75, 0.30, 0.55) & 0.0295 (0.80, 0.50, 0.50) & 0.0211 (0.75, 0.30, 0.55) & 0.019 & 4.46 \\
5 & 15 & 0.0031 (0.60, 0.70) & 0.0067 (0.70, 0.25, 0.63) & 0.0086 (0.70, 0.45, 0.63) & 0.0066 (0.70, 0.25, 0.63) & 0.035 & 5.60 \\
\addlinespace[0.5em]
10 & 5 & 0.0078 (0.65, 0.64) & 0.0219 (0.75, 0.30, 0.56) & 0.0312 (0.80, 0.55, 0.51) & 0.0207 (0.75, 0.30, 0.56) & 0.019 & 4.43 \\
10 & 10 & 0.0030 (0.60, 0.70) & 0.0066 (0.70, 0.25, 0.63) & 0.0086 (0.70, 0.45, 0.63) & 0.0065 (0.70, 0.25, 0.63) & 0.035 & 5.58 \\
10 & 15 & 0.0014 (0.55, 0.75) & 0.0025 (0.65, 0.20, 0.69) & 0.0028 (0.65, 0.35, 0.69) & 0.0024 (0.65, 0.20, 0.69) & 0.058 & 6.55 \\
\addlinespace[0.5em]
15 & 5 & 0.0030 (0.60, 0.71) & 0.0065 (0.70, 0.25, 0.64) & 0.0082 (0.75, 0.40, 0.60) & 0.0063 (0.70, 0.25, 0.64) & 0.034 & 5.52 \\
15 & 10 & 0.0014 (0.55, 0.75) & 0.0024 (0.65, 0.20, 0.69) & 0.0030 (0.70, 0.35, 0.66) & 0.0024 (0.65, 0.20, 0.69) & 0.057 & 6.51 \\
15 & 15 & 0.0007 (0.50, 0.80) & 0.0011 (0.60, 0.20, 0.74) & 0.0012 (0.65, 0.30, 0.72) & 0.0011 (0.60, 0.20, 0.74) & 0.087 & 7.31 \\
\bottomrule
\end{tabular}
}
\end{table}
\begin{table}[H]
\centering
\caption{Error bounds for $p=30$ and $p_1=20$. The values in parentheses represent the parameters used to minimize each bound. For $\mathrm{BOUND1}$, the parameters are $(v_{\min}, c_{v_{\min}})$. For $\mathrm{BOUND2}$, $\mathrm{BOUND3}$, and $\mathrm{BOUND4}$, the parameters are $(v_{\min}, c_{\min}, c_{v_{\min}})$.}
\label{tab:4}
\footnotesize
\setlength{\tabcolsep}{4pt} 
\renewcommand{\arraystretch}{1.2}
\resizebox{\textwidth}{!}{%
\begin{tabular}{rrcccccr}
\toprule
$n$ & $n_1$ & $\mathrm{BOUND1}$ & $\mathrm{BOUND2}$ & $\mathrm{BOUND3}$ & $\mathrm{BOUND4}$ & $\kappa_{\lambda}^{(2)}$ & $m$ \\
\midrule
50 & 40 & 0.0017 (0.80, 0.75) & 0.0032 (0.95, 0.20, 0.68) & 0.0040 (0.95, 0.40, 0.68) & 0.0032 (0.95, 0.20, 0.68) & 0.838 & 4.35 \\
100 & 80 & 0.0007 (0.45, 0.78) & 0.0011 (0.50, 0.20, 0.75) & 0.0012 (0.55, 0.30, 0.71) & 0.0011 (0.50, 0.20, 0.75) & 0.124 & 8.73 \\
500 & 400 & 0.0007 (0.35, 0.78) & 0.0010 (0.40, 0.15, 0.74) & 0.0011 (0.40, 0.30, 0.74) & 0.0010 (0.40, 0.15, 0.74) & 0.004 & 11.54 \\
5000 & 4000 & 0.0007 (0.30, 0.80) & 0.0011 (0.40, 0.15, 0.72) & 0.0011 (0.40, 0.30, 0.72) & 0.0011 (0.35, 0.20, 0.76) & 0.000 & 12.12 \\
\bottomrule
\end{tabular}%
}
\end{table}
\begin{table}[H]
\centering
\caption{Error bounds for $p_1/n=0.4$, $p/n_1=0.75$, and $p_1/n_1=0.5$. The values in parentheses represent the parameters used to minimize each bound. For $\mathrm{BOUND1}$, the parameters are $(v_{\min}, c_{v_{\min}})$. For $\mathrm{BOUND2}$, $\mathrm{BOUND3}$ and $\mathrm{BOUND4}$, the parameters are $(v_{\min}, c_{\min}, c_{v_{\min}})$.}
\label{tab:5}
\footnotesize
\setlength{\tabcolsep}{4pt} 
\renewcommand{\arraystretch}{1.3}
\resizebox{\textwidth}{!}{%
\begin{tabular}{rrrr r@{\hspace{2pt}}l r@{\hspace{2pt}}l r@{\hspace{2pt}}l r@{\hspace{2pt}}l cr}
\toprule
$n$ & $n_1$ & $p$ & $p_1$ & \multicolumn{2}{c}{$\mathrm{BOUND1}$} & \multicolumn{2}{c}{$\mathrm{BOUND2}$} & \multicolumn{2}{c}{$\mathrm{BOUND3}$} & \multicolumn{2}{c}{$\mathrm{BOUND4}$} & $\kappa_{\lambda}^{(2)}$ & $m$ \\
\midrule
50 & 40 & 30 & 20 & 0.0017 & (0.80, 0.75) & 0.0032 & (0.95, 0.20, 0.68) & 0.0040 & (0.95, 0.40, 0.68) & 0.0032 & (0.95, 0.20, 0.68) & 0.838 & 4.35 \\
100 & 80 & 60 & 40 & 0.0001 & (0.50, 0.86) & 0.0002 & (0.60, 0.15, 0.83) & 0.0002 & (0.60, 0.25, 0.83) & 0.0002 & (0.60, 0.15, 0.83) & 0.814 & 8.80 \\
1000 & 800 & 600 & 400 & 0.0000 & (0.10, 0.98) & 0.0000 & (0.10, 0.05, 0.98) & 0.0000 & (0.10, 0.10, 0.98) & 0.0000 & (0.10, 0.05, 0.98) & 0.791 & 88.74 \\
5000 & 4000 & 3000 & 2000 & 0.0000 & (0.05, 0.99) & 0.0000 & (0.05, 0.05, 0.99) & 0.0000 & (0.05, 0.05, 0.99) & 0.0000 & (0.05, 0.05, 0.99) & 0.789 & 444.02 \\
\bottomrule
\end{tabular}
}
\end{table}
\subsection{Cornish-Fisher expansions}
\label{subsec:Cornish-Fisher}
\begin{figure}[H]
\centering
\includegraphics[width=0.9\columnwidth]{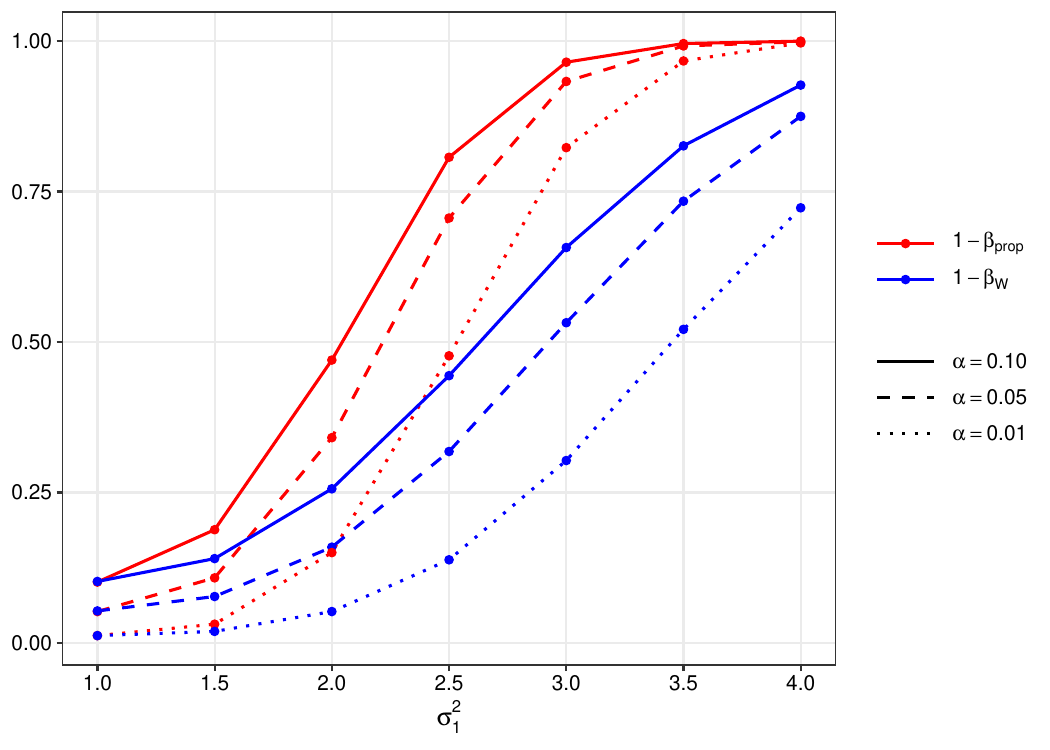}
\caption{Line graph of the empirical power $1-\beta_{\mathrm{prop}}$ and $1-\beta_{\mathrm{W}}$ for $\sigma_1^2=1.0,1.5,2.0,2.5,3.0,3.5,4.0$. The sample sizes and dimensions for $1-\beta_{\mathrm{prop}}$ are $N_1=N_2=100$ and $p_1=p_2=20$. Furthermore, those for $1-\beta_{\mathrm{W}}$ are $N_1=100$, $N_2=0$, and $p=40$. The red lines represent $1-\beta_{\mathrm{prop}}$ and blue lines represent $1-\beta_{\mathrm{W}}$. The solid, dashed, and dotted lines indicate $\alpha=0.10,0.05$, and $0.01$, respectively.}
\label{fig:3}
\end{figure}

We derive our high-dimensional Cornish-Fisher expansion for the null distribution of the LRT statistic. 
\begin{thm}\label{thm:3}
The Cornish-Fisher expansion for the null distribution of the LRT statistic up to the order $O(m^{-2})$ under high-dimensional two-step monotone incomplete data is given by
\begin{align*}
h_{\alpha}=z_{\alpha}+\dfrac{1}{m}\dfrac{1}{6}\widetilde{\kappa}_{T}^{(3)}(z_{\alpha}^2-1)+\dfrac{1}{m^2}\left\{\dfrac{1}{24}\widetilde{\kappa}_{T}^{(4)}z_{\alpha}(z_{\alpha}^2-3)+\dfrac{1}{36}(\widetilde{\kappa}_{T}^{(3)})^2z_{\alpha}(5-2z_{\alpha}^2)\right\},
\end{align*}
where $z_{\alpha}$ denotes the $1-\alpha$ quantile of the standard normal distribution and $m$ is given by Lemma \ref{lem:1}. 
\end{thm}
See \ref{app:C} for the proof of Theorem \ref{thm:3}. 
From the result of Theorem \ref{thm:3}, we numerically evaluate the empirical Type I error rates using $10^6$ Monte Carlo simulations. The settings for $\bm\mu$, $\bm\Sigma$, and $\alpha$ are identical to those in Section \ref{subsec:Type I error rates}.
\begin{align*}
\alpha_{\mathrm{prop}}&=\mathrm{Pr}\{T>h_{\alpha}\},\\
    \alpha_{\mathrm{SYS1}}&=\mathrm{Pr}\{-2\log\lambda>\ell_1(\alpha)\},\alpha_{\mathrm{SYS2}}=\mathrm{Pr}\{-2\rho\log\lambda>\ell_2(\alpha)\},
\end{align*}
where 
\begin{align*}
\ell_1(\alpha)&=\chi_f^{2}(\alpha)+\dfrac{1}{N}\dfrac{2\beta}{f}\chi_f^{2}(\alpha)+\dfrac{1}{N^2}\dfrac{2\gamma}{f(f+2)}\chi_f^{2}(\alpha)\{\chi_f^{2}(\alpha)+f+2\},\\
\ell_2(\alpha)&=\chi_f^{2}(\alpha)+\dfrac{1}{M^2}\dfrac{2\gamma^{*}}{f(f+2)}\chi_f^{2}(\alpha)\{\chi_f^{2}(\alpha)+f+2\}.
\end{align*}
$f,M,\beta,\gamma,\rho$, and $\gamma^{*}$ are given by Theorem \ref{thm:1}.
$\ell_1(\alpha)$ and $\ell_2(\alpha)$ denote the large-sample Cornish-Fisher expansions for the null distributions of the LRT statistic and modified LRT statistic up to the orders $O(N^{-2})$ and $O(M^{-2})$ under two-step monotone incomplete data given by \cite{Sato2025sphericity}. Furthermore, as a baseline for comparison with $\alpha_{\mathrm{prop}}$, let $\alpha_{\mathrm{W}}$ denote the Type I error rate in the high-dimensional complete-data subset corresponding to $N_2=0$, which consists solely of the $N_1\times p$ data matrix. Specifically, $\alpha_{\mathrm{W}}$ is obtained from the Cornish-Fisher expansion for the null distribution of the LRT statistic up to the order $O(m^{-2})$ under high-dimensional complete data, following \cite{wakaki2007error}, \cite{akita2010high}, and \cite{Ulyanov2016}.
Table \ref{tab:6} presents the empirical Type I error rates $\alpha_1$, $\alpha_{\mathrm{prop}}$, $\alpha_{\mathrm{SYS1}}$, $\alpha_{\mathrm{SYS2}}$, and $\alpha_{\mathrm{W}}$. The sample sizes and dimensions for the first four rates are $N_1=N_2=100$ and $p_1=p_2=2,10,20,40,45$. Furthermore, those for $\alpha_{\mathrm{W}}$ are $N_1=100$, $N_2=0$, and $p=4,20,40,80,90$. 
The nominal significance levels are $\alpha=0.10,0.05,0.01$. Table \ref{tab:6} indicates that $\alpha_1$, $\alpha_{\mathrm{SYS1}}$, and $\alpha_{\mathrm{SYS2}}$ deviate from $\alpha$ as $p_1$ and $p_2$ increase, whereas $\alpha_{\mathrm{prop}}$ and $\alpha_{\mathrm{W}}$ are precisely controlled. A comparison of $\alpha_{\mathrm{prop}}$ and $\alpha_{\mathrm{W}}$ indicates that $\alpha_{\mathrm{W}}$ is slightly higher in low dimensions, but this difference becomes negligible as the dimensions increase. Consequently, we also numerically evaluate the empirical power in high-dimensional settings.
\begin{table}[H]
\centering
\caption{Empirical Type I error rates $\alpha_1$, $\alpha_{\mathrm{prop}}$, $\alpha_{\mathrm{SYS1}}$, $\alpha_{\mathrm{SYS2}}$, and $\alpha_{\mathrm{W}}$. The sample sizes and dimensions for $\alpha_1$, $\alpha_{\mathrm{prop}}$, $\alpha_{\mathrm{SYS1}}$, and $\alpha_{\mathrm{SYS2}}$ are $N_1=N_2=100$ and $p_1=p_2=2,10,20,40,45$. Furthermore, those for $\alpha_{\mathrm{W}}$ are $N_1=100$, $N_2=0$, and $p=4,20,40,80,90$. The nominal significance levels are $\alpha=0.10,0.05,0.01$.}
\label{tab:6}
\renewcommand{\arraystretch}{1.2}
\footnotesize
\setlength{\tabcolsep}{16pt}
\begin{tabular}{rrrccccc}
\toprule
$p_1$ & $p_2$ & $p$ & $\alpha_1$ & $\alpha_{\mathrm{prop}}$ & $\alpha_{\mathrm{SYS1}}$ & $\alpha_{\mathrm{SYS2}}$ & $\alpha_{\mathrm{W}}$ \\
\midrule
\multicolumn{8}{c}{$\alpha=0.10$} \\
\midrule
2  & 2  & 4 & 0.111 & 0.100 & 0.100 & 0.100 & 0.102\\
10 & 10 & 20 & 0.325 & 0.103 & 0.102 & 0.101 & 0.103\\
20 & 20 & 40 & 0.973 & 0.102 & 0.134 & 0.106 & 0.102\\
40 & 40 & 80 & 1.000 & 0.101 & 0.986 & 0.513 & 0.101\\
45 & 45 & 90 & 1.000 & 0.101 & 1.000 & 0.952 & 0.101\\
\midrule
\multicolumn{8}{c}{$\alpha = 0.05$} \\
\midrule
2  & 2  & 4 & 0.057 & 0.050 & 0.050 & 0.050 & 0.057\\
10 & 10 & 20 & 0.210 & 0.054 & 0.051 & 0.050 & 0.055\\
20 & 20 & 40 & 0.943 & 0.052 & 0.070 & 0.053 & 0.053\\
40 & 40 & 80 & 1.000 & 0.052 & 0.969 & 0.377 & 0.052\\
45 & 45 & 90 & 1.000 & 0.052 & 1.000 & 0.909 & 0.052\\
\midrule
\multicolumn{8}{c}{$\alpha = 0.01$} \\
\midrule
2  & 2  & 4 & 0.012 & 0.010 & 0.010 & 0.010 & 0.016\\
10 & 10 & 20 & 0.071 & 0.013 & 0.010 & 0.010 & 0.013\\
20 & 20 & 40 & 0.827 & 0.012 & 0.016 & 0.011 & 0.012\\
40 & 40 & 80 & 1.000 & 0.011 & 0.889 & 0.169 & 0.011\\
45 & 45 & 90 & 1.000 & 0.011 & 1.000 & 0.765 & 0.011\\
\bottomrule
\end{tabular}
\end{table}
The power under high-dimensional two-step monotone incomplete data is given by $1-\beta_{\mathrm{prop}}=\mathrm{Pr}\{T>h_{\alpha}\}$,
and let $1-\beta_{\mathrm{W}}$ denote the power in the high-dimensional complete-data subset corresponding to $N_2=0$. $1-\beta_{\mathrm{W}}$ is obtained from the Cornish-Fisher expansion for the null distribution of the LRT statistic up to the order $O(m^{-2})$ under high-dimensional complete data, following \cite{wakaki2007error}, \cite{akita2010high}, and \cite{Ulyanov2016}.
Figure \ref{fig:3} displays the line graph of the empirical power $1-\beta_{\mathrm{prop}}$ and $1-\beta_{\mathrm{W}}$ when $\bm{\mu}=\bm{0},\bm{\Sigma}=\operatorname{diag}(\sigma_1^2,1,\dots,1)
$, and $\sigma_1^2=1.0,1.5,2.0,2.5,3.0,3.5,4.0$. The sample sizes and dimensions for $1-\beta_{\mathrm{prop}}$ are $N_1=N_2=100$ and $p_1=p_2=20$. Furthermore, those for $1-\beta_{\mathrm{W}}$ are $N_1=100$, $N_2=0$, and $p=40$.
The nominal significance levels are $\alpha=0.10,0.05,0.01$.
Figure \ref{fig:3} indicates that both $1-\beta_{\mathrm{prop}}$ and $1-\beta_{\mathrm{W}}$ increase with $\sigma_1^2$. Furthermore, $1-\beta_{\mathrm{prop}}$ is higher than $1-\beta_{\mathrm{W}}$ when $\sigma_1^2>1.0$. Therefore, we conclude that using partially observed two-step monotone incomplete data is preferable.
\section{Conclusions}
\label{sec:Conclusions}
This paper has derived the Edgeworth and Cornish-Fisher expansions for the null distribution of the LRT statistic for the sphericity test \eqref{eqn:a1} under high-dimensional two-step monotone incomplete data, and established the computable error bounds.
We conducted Monte Carlo simulations to compare the approximation accuracy of the Edgeworth expansion $Q_2(x)$ for $s=2$ and the large-sample asymptotic expansion for the null distribution of the LRT statistic and demonstrated that $Q_2(x)$ provides superior approximation accuracy when $p$ increases. 
Furthermore, we provide the four error bounds $\mathrm{BOUND1}$--$\mathrm{BOUND4}$.
We then evaluated these error bounds through numerical experiments. 
As $n_1$ and $p$ increase, $m$ grows and these error bounds converge to zero. 
These results are consistent with our theoretical findings and are corroborated by the evidence in Section \ref{subsec:Example}.
Finally, we show that the Type I error rate $\alpha_{\mathrm{prop}}$ in high-dimensional settings is precisely controlled.
When $n_1<p$, the theoretical framework presented in this study does not apply. Although settings where the dimension exceeds the sample size have been addressed by \cite{ledoit2002some} and \cite{Srivastava2006}, these existing approaches assume complete data. However, the case where the dimension is larger than the sample size under monotone incomplete data has not been discussed, and extending our approach to such settings can be cited as future work.
\section*{Acknowledgments}
This work was JSPS Grant-in-Aid for Early-Career Scientists Grant Number JP23K13019. The authors would like to express their deep gratitude to Professor Takashi Seo of Tokyo University of Science for his invaluable support and valuable comments throughout the course of this research.

\appendix
\makeatletter
\@addtoreset{equation}{section}
\@addtoreset{table}{section}
\makeatother
\renewcommand{\theequation}{\Alph{section}.\arabic{equation}}
\renewcommand{\thetable}{\Alph{section}.\arabic{table}}
\renewcommand{\thesection}{Appendix~\Alph{section}}
\section{Proof of Lemma \ref{lem:1}}
\label{app:A}
First, we prove the lower bound in \eqref{eqn:a4} for the high-dimensional two-step monotone incomplete data case, assuming $n>n_1>p>p_1$ and $0<\tau_1<1$. For $s\ge2$, set $H_s(x)=(-1)^s\psi^{(s-1)}(x)$, which is positive, monotonically decreasing in $x$, and convex. Then,
\begin{align*}
\kappa_{\lambda}^{(s)}=\sum_{\ell=1}^{p_1}H_s(a_{1\ell})+\tau_1^s\sum_{\ell=1}^{p_2}H_s(a_{2\ell})-(p_1+\tau_1p_2)^{s}H_s(a_3),
\end{align*}
where
\begin{align*}
a_{1\ell}=\dfrac{1}{2}(n-p_1+\ell),a_{2\ell}=\dfrac{1}{2}(n_1-p+\ell),
a_3=\dfrac{1}{2}(np_1+n_1p_2).
\end{align*}
Applying Jensen's inequality, we obtain
\begin{align*}
\sum_{\ell=1}^{p_1}H_s(a_{1\ell})>p_1H_s(\bar a_1),~\sum_{\ell=1}^{p_2}H_s(a_{2\ell})>p_2H_s(\bar a_2),
\end{align*}
where
\begin{align*}
\bar a_1=\dfrac{1}{p_1}\sum_{\ell=1}^{p_1}a_{1\ell},\bar a_2=\dfrac{1}{p_2}\sum_{\ell=1}^{p_2}a_{2\ell}.
\end{align*}
Hence, we have 
\begin{align*}
\kappa_{\lambda}^{(s)}>p_1H_s(\bar a_1)+\tau_1^sp_2H_s(\bar a_2)-(p_1+\tau_1p_2)^{s}H_s(a_3).
\end{align*}
For $s\ge2$ and $a>0$, the function $H_s(a)$ satisfies
\begin{align*}
H_s(a)=\sum_{x=0}^{\infty}\dfrac{(s-1)!}{(a+x)^s}&=\dfrac{(s-1)!}{a^s}+\sum_{x=1}^{\infty}\dfrac{(s-1)!}{(a+x)^s}\\
&<\dfrac{(s-1)!}{a^s}+\int_{0}^{\infty}\dfrac{(s-1)!}{(a+x)^s}dx=\dfrac{(s-1)!}{a^s}+\dfrac{(s-2)!}{a^{s-1}},\\
H_s(a)&>\dfrac{(s-1)!}{2a^s}+\int_{0}^{\infty}\dfrac{(s-1)!}{(a+x)^s}dx=\dfrac{(s-1)!}{2a^s}+\dfrac{(s-2)!}{a^{s-1}}.
\end{align*}
Hence, we have
\begin{align}
\kappa_{\lambda}^{(s)}>&~p_1\biggl\{\dfrac{(s-1)!}{2{\bar a_1}^s}+\dfrac{(s-2)!}{{\bar a_1}^{s-1}}\biggr\}+\tau_1^sp_2\biggl\{\dfrac{(s-1)!}{2{\bar a_2}^s}+\dfrac{(s-2)!}{{\bar a_2}^{s-1}}\biggr\}\\
&-(p_1+\tau_1p_2)^{s}\biggl\{\dfrac{(s-1)!}{a_3^s}+\dfrac{(s-2)!}{a_3^{s-1}}\biggr\}\\
=&~\dfrac{(s-1)!}{2}D_1+(s-2)!D_2,
\end{align}
where
\begin{align*}
D_1&=p_1{\bar a_1}^{-s}+\tau_1^sp_2{\bar a_2}^{-s}-2(p_1+\tau_1p_2)^{s}a_3^{-s},\\
D_2&=p_1{\bar a_1}^{1-s}+\tau_1^sp_2{\bar a_2}^{1-s}-(p_1+\tau_1p_2)^{s}a_3^{1-s}.
\end{align*}
We can write $\bar a_1$, $\bar a_2$, and $a_3$ as
\begin{align*}
\bar a_1=\dfrac{N}{2}(1-c_1),\bar a_2=\dfrac{N_1}{2}(1-c_2),a_3=\dfrac{N}{2}(p_1+\tau_1p_2)(1-c_3),
\end{align*}
respectively, where
\begin{align*}
c_1=\dfrac{p_1+1}{2N},c_2=\dfrac{p+p_1+1}{2N_1},c_3=\dfrac{p}{N(p_1+\tau_1p_2)}.
\end{align*}
Hence, $D_1$ can be written in the form
\begin{align*}
D_1&=\left(\dfrac{2}{N}\right)^{s}\{p_1g_1(c_1)+p_2g_1(c_2)-2g_1(c_3)\}\\
&=\left(\dfrac{2}{N}\right)^{s}\Bigl\{p_1(g_1(c_1)-g_1(c_3))+p_2(g_1(c_2)-g_1(c_3))
+(p-2)g_1(c_3)\Bigr\},
\end{align*}
where $g_1(c_i)=(1-c_i)^{-s}$ for $c_i\in(0,1)$ and $i=1,2,3$.
Let $c_{max}=\max\{c_1,c_2,c_3\}$, then $c_{max}=c_2$. By the mean value theorem, we have
\begin{align*}
g_1(c_j)-g_1(c_3)>s(1-c_2)^{-s-1}(c_j-c_3),
\end{align*}
for $j=1,2$. This inequality yields the lower bound
\begin{align}
\begin{split}\label{eqn:A.1}
D_1>&\left(\dfrac{2}{N}\right)^{s}\Biggl[s(1-c_2)^{-s-1}\Bigl\{\dfrac{1}{2N}p(p-1)+(1-\tau_1)p_2(c_2-c_3)\Bigr\}\\
&+(p-2)(1-c_3)^{-s}\Biggr].
\end{split}
\end{align}
Applying the same arguments to $D_2$ yields
\begin{align}
\begin{split}\label{eqn:A.2}
D_2&=\left(\dfrac{2}{N}\right)^{s-1}\{p_1(g_2(c_1)-g_2(c_3))+\tau_1p_2(g_2(c_2)-g_2(c_3))\}\\
&>\dfrac{1}{4}\left(\dfrac{2}{N}\right)^{s}(s-1)(1-c_2)^{-s}p(p-1),
\end{split}
\end{align}
where $g_2(c_i)=(1-c_i)^{-(s-1)}$ for $c_i\in(0,1)$ and $i=1,2,3$.
From \eqref{eqn:A.1} and \eqref{eqn:A.2}, $D_1$ and $D_2$ are positive. Therefore, we obtain  $\kappa_{\lambda}^{(s)}>0$ for $s\ge 2$.
In the complete data case where $N_2=0$, it follows from \cite{wakaki2007error} that $\kappa_{\lambda}^{(s)}$ is positive for $s\ge 2$. Therefore, the lower bound in \eqref{eqn:a4} is proved.

We next prove the upper bound in \eqref{eqn:a4}. We have
\begin{align*}
H_s(a)<\int_{-\tfrac{1}{2}}^{\infty}\dfrac{(s-1)!}{(a+x)^s}dx=\dfrac{(s-2)!}{(a-\tfrac{1}{2})^{s-1}},
\end{align*}
for $a>1/2$, and
\begin{align*}
H_s(a)>\dfrac{(s-1)!}{2a^s}+\dfrac{(s-2)!}{a^{s-1}},
\end{align*}
for $a>0$.
Therefore, the following inequalities hold for $s\ge3$.
\begin{align}
\begin{split}\label{eqn:A.3}
&(-1)^s\sum_{\ell=1}^{p_1}\psi^{(s-1)}\biggl(\dfrac{1}{2}(n-p_1+\ell)\biggr)\\
<&~\sum_{\ell=1}^{p_1}\dfrac{(s-2)!2^{s-1}}{(n-p_1-1+\ell)^{s-1}}\\
<&\int_{\tfrac{1}{2}}^{p_1+\tfrac{1}{2}}\dfrac{(s-2)!2^{s-1}}{(n-p_1-1+x)^{s-1}}dx\\
=&~\dfrac{(s-3)!2^{s-1}}{(n-p_1-\tfrac{1}{2})^{s-2}}\Biggl\{1-\Bigl(\dfrac{n-p_1-\tfrac{1}{2}}{n-\tfrac{1}{2}}\Bigr)^{s-2}\Biggr\},
\end{split}
\end{align}
\begin{align}
\begin{split}\label{eqn:A.4}
&(-1)^s\tau_1^s\sum_{\ell=1}^{p_2}\psi^{(s-1)}\biggl(\dfrac{1}{2}(n_1-p+\ell)\biggr)\\
<&~\tau_1^s\sum_{\ell=1}^{p_2}\dfrac{(s-2)!2^{s-1}}{(n_1-p-1+\ell)^{s-1}}\\
<&~\tau_1^s\int_{\tfrac{1}{2}}^{p_2+\tfrac{1}{2}}\dfrac{(s-2)!2^{s-1}}{(n_1-p-1+x)^{s-1}}dx\\
=&~\dfrac{\tau_1^s(s-3)!2^{s-1}}{(n_1-p-\tfrac{1}{2})^{s-2}}\Biggl\{1-\Bigl(\dfrac{n_1-p-\tfrac{1}{2}}{n_1-p_1-\tfrac{1}{2}}\Bigr)^{s-2}\Biggr\},
\end{split}
\end{align}
\begin{align}
\begin{split}\label{eqn:A.5}
&(-1)^{s-1}(p_1+\tau_1p_2)^{s}\psi^{(s-1)}\biggl(\dfrac{1}{2}(np_1+n_1p_2)\biggr)\\
<&-(p_1+\tau_1p_2)^{s}\Biggl[\dfrac{(s-1)!}{2\{\tfrac{1}{2}(np_1+n_1p_2)\}^s}+\dfrac{(s-2)!}{\{\tfrac{1}{2}(np_1+n_1p_2)\}^{s-1}}\Biggr]\\
=&-\dfrac{(p_1+\tau_1p_2)^{s}(s-1)!2^{s-1}}{(np_1+n_1p_2)^s}-\dfrac{(p_1+\tau_1p_2)^{s}(s-2)!2^{s-1}}{(np_1+n_1p_2)^{s-1}}.
\end{split}
\end{align}
Moreover, the upper bound in \eqref{eqn:A.5} satisfies the following inequality.
\begin{align*}
&-\dfrac{(p_1+\tau_1p_2)^{s}(s-1)!2^{s-1}}{(np_1+n_1p_2)^s}
-\dfrac{(p_1+\tau_1p_2)^{s}(s-2)!2^{s-1}}{(np_1+n_1p_2)^{s-1}}\\
<&-\dfrac{(p_1+\tau_1p_2)^{s}(s-1)!2^{s-1}}{\{n(p_1+\tau_1p_2)\}^s}-\dfrac{(p_1+\tau_1p_2)^{s}(s-2)!2^{s-1}}{\{n(p_1+\tau_1p_2)\}^{s-1}}\\
=&-\dfrac{(s-1)!2^{s-1}}{n^s}-\dfrac{(p_1+\tau_1p_2)(s-2)!2^{s-1}}{n^{s-1}}.
\end{align*}
Therefore,
\begin{align}
\begin{split}\label{eqn:A.6}
&(-1)^{s-1}(p_1+\tau_1p_2)^{s}\psi^{(s-1)}\biggl(\dfrac{1}{2}(np_1+n_1p_2)\biggr)\\
<&-\dfrac{(s-1)!2^{s-1}}{n^s}-\dfrac{(p_1+\tau_1p_2)(s-2)!2^{s-1}}{n^{s-1}}.
\end{split}
\end{align}
From \eqref{eqn:A.3}, \eqref{eqn:A.4}, and \eqref{eqn:A.6},
\begin{align*}
\kappa_{\lambda}^{(s)}<&-\dfrac{(s-1)!2^{s-1}}{n^s}-\dfrac{(p_1+\tau_1p_2)(s-2)!2^{s-1}}{n^{s-1}}\\
&+\dfrac{(s-3)!2^{s-1}}{(n-p_1-\tfrac{1}{2})^{s-2}}\Biggl\{1-\Bigl(\dfrac{n-p_1-\tfrac{1}{2}}{n-\tfrac{1}{2}}\Bigr)^{s-2}\Biggr\}\\
&+\dfrac{\tau_1^s(s-3)!2^{s-1}}{(n_1-p-\tfrac{1}{2})^{s-2}}\Biggl\{1-\Bigl(\dfrac{n_1-p-\tfrac{1}{2}}{n_1-p_1-\tfrac{1}{2}}\Bigr)^{s-2}\Biggr\},
\end{align*}
for $s\ge3$. Hence, for $s\ge3$,
\begin{align*}
\dfrac{\widetilde{\kappa}_{T}^{(s)}}{s!}=&~\dfrac{\kappa_{\lambda}^{(s)}}{(\kappa_{\lambda}^{(2)})^{\tfrac{s}{2}}}\dfrac{1}{s!}\\
<&-\dfrac{2^{s-1}}{sn^s(\kappa_{\lambda}^{(2)})^{\tfrac{s}{2}}}-\dfrac{(p_1+\tau_1p_2)2^{s-1}}{(s-1)sn^{s-1}(\kappa_{\lambda}^{(2)})^{\tfrac{s}{2}}}\\
&+\dfrac{2^{s-1}}{(s-2)(s-1)s(n-p_1-\tfrac{1}{2})^{s-2}(\kappa_{\lambda}^{(2)})^{\tfrac{s}{2}}}\Biggl\{1-\Bigl(\dfrac{n-p_1-\tfrac{1}{2}}{n-\tfrac{1}{2}}\Bigr)^{s-2}\Biggr\}\\
&+\dfrac{\tau_1^s2^{s-1}}{(s-2)(s-1)s(n_1-p-\tfrac{1}{2})^{s-2}(\kappa_{\lambda}^{(2)})^{\tfrac{s}{2}}}\Biggl\{1-\Bigl(\dfrac{n_1-p-\tfrac{1}{2}}{n_1-p_1-\tfrac{1}{2}}\Bigr)^{s-2}\Biggr\}\\
=&~\dfrac{2^{s-2}}{(n_1-p-\tfrac{1}{2})^{s-2}(\kappa_{\lambda}^{(2)})^{\tfrac{s-2}{2}}}
\Biggl[\dfrac{2}{\kappa_{\lambda}^{(2)}(s-2)(s-1)s}\\
&\times\Biggl\{\biggl(\dfrac{n_1-p-\tfrac{1}{2}}{n-p_1-\tfrac{1}{2}}\biggr)^{s-2}-\biggl(\dfrac{n_1-p-\tfrac{1}{2}}{n-\tfrac{1}{2}}\biggr)^{s-2}\Biggr\}\\
&+\dfrac{2\tau_1^2}{\kappa_{\lambda}^{(2)}(s-2)(s-1)s}\Biggl\{\tau_1^{s-2}-\biggl(\tau_1\dfrac{n_1-p-\tfrac{1}{2}}{n_1-p_1-\tfrac{1}{2}}\biggr)^{s-2}\Biggr\}\\
&-\Biggl\{\dfrac{2}{\kappa_{\lambda}^{(2)}sn^2}+\dfrac{2(p_1+\tau_1p_2)}{\kappa_{\lambda}^{(2)}(s-1)sn}\Biggr\}\biggl(\dfrac{n_1-p-\tfrac{1}{2}}{n}\biggr)^{s-2}\Biggr]=m^{-(s-2)}b_{s-3},
\end{align*}
where $m$ and $b_s$ are defined as in Lemma \ref{lem:1}.
Hence, the upper bound in \eqref{eqn:a4} has been established, which completes the proof.
\section{Proof of \texorpdfstring{$\gamma_{k,j}=O(m^{-(j+k)})$}{gamma kj = O(m)}} 
\label{app:B}
By Lemma \ref{lem:1}, we have
\begin{align*}
\dfrac{|\widetilde{\kappa}_{T}^{(s_1+3)}|\cdots|\widetilde{\kappa}_{T}^{(s_k+3)}|}{(s_1+3)!\cdots(s_k+3)!}&<m^{-(s_1+1)}\cdots m^{-(s_k+1)}b_{s_1}\cdots b_{s_k}=m^{-(s_1+\cdots+s_k+k)}b_{s_1}\cdots b_{s_k}.
\end{align*}
Consequently,
\begin{align*}
|\gamma_{k,j}|\leq\sum_{s_1+\cdots+s_k=j}\dfrac{|\widetilde{\kappa}_{T}^{(s_1+3)}|\cdots|\widetilde{\kappa}_{T}^{(s_k+3)}|}{(s_1+3)!\cdots(s_k+3)!}
<m^{-(j+k)}\sum_{s_1+\cdots+s_k=j}b_{s_1}\cdots b_{s_k}.
\end{align*}
By the definition of big-$O$ notation, we obtain
\begin{align*}
\gamma_{k,j}=O(m^{-(j+k)})\qquad(m\rightarrow\infty).
\end{align*}
\section{Proof of Theorem \ref{thm:3}}
\label{app:C}
We express the Edgeworth expansion $Q_2(x)$ and Cornish-Fisher expansion $h_\alpha$ up to the order $O(m^{-2})$ as
\begin{align*}
Q_2(x)=\Phi(x)-\phi(x)\left(\dfrac{1}{m}L_1(x)+\dfrac{1}{m^2}L_2(x)\right)
\end{align*}
and
\begin{align*}
h_\alpha=z_\alpha+\dfrac{1}{m}\varepsilon_1+\dfrac{1}{m^2}\varepsilon_2,
\end{align*}
respectively, where
\begin{align*}
L_1(x)&=\dfrac{1}{6}\widetilde{\kappa}_{T}^{(3)}(x^2-1),
L_2(x)=\dfrac{1}{24}\widetilde{\kappa}_{T}^{(4)}(x^3-3x)+\dfrac{1}{72}(\widetilde{\kappa}_{T}^{(3)})^2(x^5-10x^3+15x).
\end{align*}
Based on the generalized Cornish-Fisher expansion derived by \cite{Ulyanov2016}, we obtain
\begin{align}
\varepsilon_1&=L_1(z_\alpha)=\dfrac{1}{6}\widetilde{\kappa}_{T}^{(3)}(z_\alpha^2-1),\\
\varepsilon_2&=-\dfrac{1}{2}\dfrac{\phi'(z_\alpha)}{\phi(z_\alpha)}L_1(z_\alpha)+L_2(z_\alpha)+L_1'(z_\alpha)L_1(z_\alpha)\\
&=\dfrac{1}{24}\widetilde{\kappa}_{T}^{(4)}z_\alpha(z_\alpha^2-3)+\dfrac{1}{36}(\widetilde{\kappa}_{T}^{(3)})^2z_\alpha(5-2z_\alpha^2).
\end{align}
This yields Theorem \ref{thm:3}.
\section{Comparison of approximation accuracy under various parameters}
\label{app:D}
Tables \ref{tab:D.1}--\ref{tab:D.4} report the biases
$B_{\mathrm{prop}}$ and $B_{\mathrm{SYS}}$ under the following parameter settings: 
$(i)~N_1=50,100,200$, $N_1/N_2=0.5,1,2$, and $p_1=p_2=2$; and $(ii)~N_1=50,100,200$, $N_1/N_2=0.5,1,2$, $p/N_1=0.2,0.4,0.8$, and $p_1/p_2=0.25,1,4$.
Tables \ref{tab:D.2}--\ref{tab:D.4} present the results for $\alpha=0.10,0.05$, and $0.01$, respectively. Table \ref{tab:D.1} shows that the large-sample asymptotic expansion for the null distribution of the likelihood ratio test (LRT) statistic given by \eqref{eqn:a2} provides a slightly more accurate approximation than the Edgeworth expansion $Q_2(x)$, whereas Tables \ref{tab:D.2}--\ref{tab:D.4} indicate that $Q_2(x)$ is clearly superior.
\begin{table}[H]
\centering
\caption{Biases $B_{\mathrm{prop}}$ and $B_{\mathrm{SYS}}$ for $(i)~N_1=50,100,200$, $N_1/N_2=0.5,1,2$, $p_1=p_2=2$, and $\alpha=0.10,0.05,0.01$. In each row, whichever of $B_{\mathrm{prop}}$ and $B_{\mathrm{SYS}}$ is closer to zero is shown in bold.}
\label{tab:D.1}
\renewcommand{\arraystretch}{1.1}
\setlength{\tabcolsep}{16pt}
\begin{tabular}{rrrr}
\toprule
$N_1$ & $N_2$ & $B_{\mathrm{prop}}$ & $B_{\mathrm{SYS}}$\\
\midrule
\multicolumn{4}{c}{$\alpha=0.10$} \\
\midrule
50  & 50  & 0.003 & \textbf{0.001} \\
100 & 100 & 0.003 & \textbf{0.000} \\
200 & 200 & 0.003 & \textbf{0.000} \\
50  & 100 & 0.003 & \textbf{0.001} \\
100 & 200 & 0.003 & \textbf{0.000} \\
200 & 400 & 0.003 & \textbf{0.000} \\
50  & 25  & 0.004 & \textbf{0.002} \\
100 & 50  & 0.003 & \textbf{0.000} \\
200 & 100 & 0.003 & \textbf{0.000} \\
\midrule
\multicolumn{4}{c}{$\alpha = 0.05$} \\
\midrule
50  & 50  & 0.002 & \textbf{0.001} \\
100 & 100 & 0.002 & \textbf{0.000} \\
200 & 200 & 0.001 & \textbf{0.000} \\
50  & 100 & 0.002 & \textbf{0.001} \\
100 & 200 & 0.001 & \textbf{0.000} \\
200 & 400 & 0.001 & \textbf{0.000} \\
50  & 25  & 0.003 & \textbf{0.002} \\
100 & 50  & 0.001 & \textbf{0.000} \\
200 & 100 & 0.001 & \textbf{0.000} \\
\midrule
\multicolumn{4}{c}{$\alpha = 0.01$} \\
\midrule
50  & 50  & $-0.004$ & \textbf{0.000} \\
100 & 100 & $-0.003$ & \textbf{0.000} \\
200 & 200 & $-0.003$ & \textbf{0.000} \\
50  & 100 & $-0.003$ & \textbf{0.000} \\
100 & 200 & $-0.003$ & \textbf{0.000} \\
200 & 400 & $-0.003$ & \textbf{0.000} \\
50  & 25  & $-0.003$ & \textbf{0.001} \\
100 & 50  & $-0.003$ & \textbf{0.000} \\
200 & 100 & $-0.003$ & \textbf{0.000} \\
\bottomrule
\end{tabular}
\end{table}
\begin{table}[H]
\centering
\caption{Biases $B_{\mathrm{prop}}$ and $B_{\mathrm{SYS}}$ for $(ii)~N_1=50,100,200$, $N_1/N_2=0.5,1,2$, $p/N_1=0.2,0.4,0.8$, $p_1/p_2=0.25,1,4$, and $\alpha=0.10$. In each row, whichever of $B_{\mathrm{prop}}$ and $B_{\mathrm{SYS}}$ is closer to zero is shown in bold.}
\label{tab:D.2}
\renewcommand{\arraystretch}{0.95}
\setlength{\tabcolsep}{20pt}
\begin{tabular}{rrrrrr}
\toprule
$N_1$ & $N_2$ & $p_1$ & $p_2$ & $B_{\mathrm{prop}}$ & $B_{\mathrm{SYS}}$ \\
\midrule
50  & 50 & 2 & 8 & \textbf{0.001} & 0.029 \\
50  & 50 & 5 & 5 & \textbf{0.000} & 0.025 \\
50  & 50 & 8 & 2 & \textbf{0.000} & 0.015 \\
50  & 50 & 4 & 16 & \textbf{0.000} & 0.272 \\
50  & 50 & 10 & 10 & \textbf{0.000} & 0.258 \\
50  & 50 & 16 & 4 & \textbf{0.000} & 0.195 \\
50  & 50 & 8 & 32 & \textbf{0.000} & $-0.596$ \\
50  & 50 & 20 & 20 & \textbf{0.000} & $-0.506$ \\
50  & 50 & 32 & 8 & \textbf{0.000} & $-0.183$ \\
100  & 100 & 4 & 16 & \textbf{0.000} & 0.083 \\
100  & 100 & 10 & 10 & \textbf{0.000} & 0.073 \\
100  & 100 & 16 & 4 & \textbf{0.000} & 0.046 \\
100  & 100 & 8 & 32 & \textbf{0.000} & 0.247 \\
100  & 100 & 20 & 20 & \textbf{0.000} & 0.276 \\
100  & 100 & 32 & 8 & \textbf{0.000} & 0.333 \\
100  & 100 & 16 & 64 & \textbf{0.000} & $-1.898$ \\
100  & 100 & 40 & 40 & \textbf{0.000} & $-1.737$ \\
100  & 100 & 64 & 16 & \textbf{0.000} & $-1.137$ \\
200  & 200 & 8 & 32 & \textbf{0.001} & 0.229 \\
200  & 200 & 20 & 20 & \textbf{0.001} & 0.211 \\
200  & 200 & 32 & 8 & \textbf{0.000} & 0.147 \\
200  & 200 & 16 & 64 & \textbf{0.000} & $-0.294$  \\
200  & 200 & 40 & 40 & \textbf{0.000} & $-0.223$ \\
200  & 200 & 64 & 16 & \textbf{0.000} & 0.023 \\
200  & 200 & 32 & 128 & \textbf{0.000} & $-4.509$ \\
200  & 200 & 80 & 80 & \textbf{0.000} & $-4.204$ \\
200  & 200 & 128 & 32 & \textbf{0.000} & $-3.047$ \\
50  & 100 & 2 & 8 & \textbf{0.000} & 0.029 \\
50  & 100 & 5 & 5 & \textbf{0.001} & 0.024 \\
50  & 100 & 8 & 2 & \textbf{0.001} & 0.013 \\
50  & 100 & 4 & 16 & \textbf{0.000} & 0.271 \\
50  & 100 & 10 & 10 & \textbf{0.000} & 0.254 \\
50  & 100 & 16 & 4 & \textbf{0.000} & 0.167 \\
50  & 100 & 8 & 32 & \textbf{0.000} & $-0.594$ \\
50  & 100 & 20 & 20 & \textbf{0.000} & $-0.477$ \\
50  & 100 & 32 & 8 & \textbf{0.000} & $-0.064$ \\
100  & 200 & 4 & 16 & \textbf{0.000} & 0.082 \\
100  & 200 & 10 & 10 & \textbf{0.000} & 0.070 \\
100  & 200 & 16 & 4 & $\mathbf{-0.001}$ & 0.036 \\
100  & 200 & 8 & 32 & \textbf{0.000} & 0.248 \\
100  & 200 & 20 & 20 & \textbf{0.000} & 0.285 \\
100  & 200 & 32 & 8 & \textbf{0.000} & 0.333 \\
100  & 200 & 16 & 64 & \textbf{0.000} & $-1.895$ \\
100  & 200 & 40 & 40 & \textbf{0.000} & $-1.686$ \\
100  & 200 & 64 & 16 & \textbf{0.000} & $-0.917$ \\
\bottomrule
\end{tabular}
\end{table}
\begin{table}[H]
\centering
\ContinuedFloat
\caption{(Continued)}
\setlength{\tabcolsep}{20pt}
\begin{tabular}{rrrrrr}
\toprule
$N_1$ & $N_2$ & $p_1$ & $p_2$ & $B_{\mathrm{prop}}$ & $B_{\mathrm{SYS}}$ \\
\midrule
200  & 400 & 8 & 32 & $\mathbf{-0.001}$ & 0.229 \\
200  & 400 & 20 & 20 & \textbf{0.000} & 0.204 \\
200  & 400 & 32 & 8 & \textbf{0.000} & 0.119 \\
200  & 400 & 16 & 64 & \textbf{0.000} & $-0.292$ \\
200  & 400 & 40 & 40 & \textbf{0.000} & $-0.200$ \\
200  & 400 & 64 & 16 & \textbf{0.000} & 0.120 \\
200  & 400 & 32 & 128 & \textbf{0.000} & $-4.503$ \\
200  & 400 & 80 & 80 & \textbf{0.000} & $-4.108$ \\
200  & 400 & 128 & 32 & \textbf{0.000} & $-2.625$ \\
50  & 25 & 2 & 8 & \textbf{0.000} & 0.029 \\
50  & 25 & 5 & 5 & \textbf{0.000} & 0.026 \\
50  & 25 & 8 & 2 & \textbf{0.000} & 0.019 \\
50  & 25 & 4 & 16 & \textbf{0.000} & 0.271 \\
50  & 25 & 10 & 10 & \textbf{0.000} & 0.263 \\
50  & 25 & 16 & 4 & \textbf{0.001} & 0.224 \\
50  & 25 & 8 & 32 & \textbf{0.000} & $-0.598$ \\
50  & 25 & 20 & 20 & \textbf{0.000} & $-0.536$ \\
50  & 25 & 32 & 8 & \textbf{0.000} & $-0.312$ \\
100  & 50 & 4 & 16 & $\mathbf{-0.001}$ & 0.082 \\
100  & 50 & 10 & 10 & \textbf{0.000} & 0.077 \\
100  & 50 & 16 & 4 & \textbf{0.000} & 0.058 \\
100  & 50 & 8 & 32 & \textbf{0.000} & 0.246 \\
100  & 50 & 20 & 20 & \textbf{0.000} & 0.267 \\
100  & 50 & 32 & 8 & \textbf{0.000} & 0.318 \\
100  & 50 & 16 & 64 & \textbf{0.000} & $-1.902$ \\
100  & 50 & 40 & 40 & \textbf{0.000} & $-1.791$ \\
100  & 50 & 64 & 16 & \textbf{0.000} & $-1.375$ \\
200  & 100 & 8 & 32 & \textbf{0.000} & 0.229 \\
200  & 100 & 20 & 20 & $\mathbf{-0.001}$ & 0.217 \\
200  & 100 & 32 & 8 & \textbf{0.000} & 0.175 \\
200  & 100 & 16 & 64 & \textbf{0.000} & $-0.296$ \\
200  & 100 & 40 & 40 & \textbf{0.000} & $-0.248$ \\
200  & 100 & 64 & 16 & \textbf{0.000} & $-0.080$ \\
200  & 100 & 32 & 128 & \textbf{0.000} & $-4.515$ \\
200  & 100 & 80 & 80 & \textbf{0.000} & $-4.307$ \\
200  & 100 & 128 & 32 & \textbf{0.000} & $-3.505$ \\
\bottomrule
\end{tabular}
\end{table}
\begin{table}[H]
\centering
\caption{Biases $B_{\mathrm{prop}}$ and $B_{\mathrm{SYS}}$ for $(ii)~N_1=50,100,200$, $N_1/N_2=0.5,1,2$, $p/N_1=0.2,0.4,0.8$, $p_1/p_2=0.25,1,4$, and $\alpha=0.05$. In each row, whichever of $B_{\mathrm{prop}}$ and $B_{\mathrm{SYS}}$ is closer to zero is shown in bold.}
\label{tab:D.3}
\renewcommand{\arraystretch}{0.95}
\setlength{\tabcolsep}{20pt}
\begin{tabular}{rrrrrr}
\toprule
$N_1$ & $N_2$ & $p_1$ & $p_2$ & $B_{\mathrm{prop}}$ & $B_{\mathrm{SYS}}$ \\
\midrule
50  & 50 & 2 & 8 & \textbf{0.001} & 0.027 \\
50  & 50 & 5 & 5 & \textbf{0.000} & 0.023 \\
50  & 50 & 8 & 2 & \textbf{0.000} & 0.014 \\
50  & 50 & 4 & 16 & \textbf{0.000} & 0.344 \\
50  & 50 & 10 & 10 & \textbf{0.000} & 0.315 \\
50  & 50 & 16 & 4 & $\mathbf{-0.001}$ & 0.212 \\
50  & 50 & 8 & 32 & \textbf{0.000} & 0.063 \\
50  & 50 & 20 & 20 & \textbf{0.000} & 0.117 \\
50  & 50 & 32 & 8 & \textbf{0.000} & 0.308 \\
100  & 100 & 4 & 16 & \textbf{0.000} & 0.080 \\
100  & 100 & 10 & 10 & \textbf{0.000} & 0.070 \\
100  & 100 & 16 & 4 & \textbf{0.000} & 0.042 \\
100  & 100 & 8 & 32 & \textbf{0.000} & 0.535 \\
100  & 100 & 20 & 20 & \textbf{0.000} & 0.540 \\
100  & 100 & 32 & 8 & \textbf{0.000} & 0.497 \\
100  & 100 & 16 & 64 & \textbf{0.000} & $-0.702$ \\
100  & 100 & 40 & 40 & \textbf{0.000} & $-0.606$ \\
100  & 100 & 64 & 16 & \textbf{0.000} & $-0.252$ \\
200  & 200 & 8 & 32 & \textbf{0.000} & 0.258 \\
200  & 200 & 20 & 20 & \textbf{0.000} & 0.231 \\
200  & 200 & 32 & 8 & \textbf{0.000} & 0.148 \\
200  & 200 & 16 & 64 & \textbf{0.000} & 0.246 \\
200  & 200 & 40 & 40 & \textbf{0.000} & 0.287 \\
200  & 200 & 64 & 16 & \textbf{0.000} & 0.432 \\
200  & 200 & 32 & 128 & \textbf{0.000} & $-2.236$ \\
200  & 200 & 80 & 80 & \textbf{0.000} & $-2.056$ \\
200  & 200 & 128 & 32 & \textbf{0.000} & $-1.375$ \\
50  & 100 & 2 & 8 & \textbf{0.001} & 0.027 \\
50  & 100 & 5 & 5 & \textbf{0.001} & 0.022 \\
50  & 100 & 8 & 2 & \textbf{0.000} & 0.011 \\
50  & 100 & 4 & 16 & \textbf{0.000} & 0.343 \\
50  & 100 & 10 & 10 & \textbf{0.001} & 0.307 \\
50  & 100 & 16 & 4 & \textbf{0.000} & 0.175 \\
50  & 100 & 8 & 32 & \textbf{0.000} & 0.065 \\
50  & 100 & 20 & 20 & \textbf{0.000} & 0.134 \\
50  & 100 & 32 & 8 & \textbf{0.000} & 0.378 \\
100  & 200 & 4 & 16 & \textbf{0.000} & 0.079 \\
100  & 200 & 10 & 10 & \textbf{0.001} & 0.067 \\
100  & 200 & 16 & 4 & $\mathbf{-0.001}$ & 0.032 \\
100  & 200 & 8 & 32 & \textbf{0.000} & 0.535 \\
100  & 200 & 20 & 20 & \textbf{0.000} & 0.541 \\
100  & 200 & 32 & 8 & \textbf{0.001} & 0.454 \\
100  & 200 & 16 & 64 & \textbf{0.000} & $-0.700$ \\
100  & 200 & 40 & 40 & \textbf{0.000} & $-0.576$ \\
100  & 200 & 64 & 16 & \textbf{0.000} & $-0.122$ \\
\bottomrule
\end{tabular}
\end{table}
\begin{table}[H]
\centering
\ContinuedFloat
\caption{(Continued)}
\setlength{\tabcolsep}{20pt}
\begin{tabular}{rrrrrr}
\toprule
$N_1$ & $N_2$ & $p_1$ & $p_2$ & $B_{\mathrm{prop}}$ & $B_{\mathrm{SYS}}$ \\
\midrule
200  & 400 & 8 & 32 & \textbf{0.000} & 0.257 \\
200  & 400 & 20 & 20 & \textbf{0.000} & 0.222 \\
200  & 400 & 32 & 8 & \textbf{0.000} & 0.117 \\
200  & 400 & 16 & 64 & \textbf{0.000} & 0.247 \\
200  & 400 & 40 & 40 & \textbf{0.000} & 0.301 \\
200  & 400 & 64 & 16 & \textbf{0.000} & 0.487 \\
200  & 400 & 32 & 128 & \textbf{0.000} & $-2.232$ \\
200  & 400 & 80 & 80 & \textbf{0.000} & $-1.999$ \\
200  & 400 & 128 & 32 & \textbf{0.000} & $-1.126$ \\
50  & 25 & 2 & 8 & \textbf{0.000} & 0.026 \\
50  & 25 & 5 & 5 & \textbf{0.000} & 0.024 \\
50  & 25 & 8 & 2 & \textbf{0.000} & 0.018 \\
50  & 25 & 4 & 16 & $\mathbf{-0.001}$ & 0.344 \\
50  & 25 & 10 & 10 & \textbf{0.000} & 0.325 \\
50  & 25 & 16 & 4 & \textbf{0.001} & 0.256 \\
50  & 25 & 8 & 32 & \textbf{0.000} & 0.062 \\
50  & 25 & 20 & 20 & \textbf{0.000} & 0.099 \\
50  & 25 & 32 & 8 & \textbf{0.000} & 0.232 \\
100  & 50 & 4 & 16 & \textbf{0.000} & 0.079 \\
100  & 50 & 10 & 10 & \textbf{0.000} & 0.073 \\
100  & 50 & 16 & 4 & \textbf{0.000} & 0.053 \\
100  & 50 & 8 & 32 & \textbf{0.000} & 0.535 \\
100  & 50 & 20 & 20 & \textbf{0.000} & 0.539 \\
100  & 50 & 32 & 8 & \textbf{0.000} & 0.527 \\
100  & 50 & 16 & 64 & \textbf{0.000} & $-0.704$ \\
100  & 50 & 40 & 40 & \textbf{0.000} & $-0.639$ \\
100  & 50 & 64 & 16 & \textbf{0.000} & $-0.393$ \\
200  & 100 & 8 & 32 & $\mathbf{-0.001}$ & 0.258 \\
200  & 100 & 20 & 20 & $\mathbf{-0.001}$ & 0.240 \\
200  & 100 & 32 & 8 & \textbf{0.000} & 0.183 \\
200  & 100 & 16 & 64 & \textbf{0.000} & 0.245 \\
200  & 100 & 40 & 40 & \textbf{0.000} & 0.273 \\
200  & 100 & 64 & 16 & \textbf{0.000} & 0.372 \\
200  & 100 & 32 & 128 & \textbf{0.000} & $-2.239$ \\
200  & 100 & 80 & 80 & \textbf{0.000} & $-2.116$ \\
200  & 100 & 128 & 32 & \textbf{0.000} & $-1.644$ \\
\bottomrule
\end{tabular}
\end{table}
\begin{table}[H]
\centering
\caption{Biases $B_{\mathrm{prop}}$ and $B_{\mathrm{SYS}}$ for $(ii)~N_1=50,100,200$, $N_1/N_2=0.5,1,2$, $p/N_1=0.2,0.4,0.8$, $p_1/p_2=0.25,1,4$, and $\alpha=0.01$. In each row, whichever of $B_{\mathrm{prop}}$ and $B_{\mathrm{SYS}}$ is closer to zero is shown in bold.}
\label{tab:D.4}
\renewcommand{\arraystretch}{0.95}
\setlength{\tabcolsep}{20pt}
\begin{tabular}{rrrrrr}
\toprule
$N_1$ & $N_2$ & $p_1$ & $p_2$ & $B_{\mathrm{prop}}$ & $B_{\mathrm{SYS}}$ \\
\midrule
50  & 50 & 2 & 8 & \textbf{0.000} & 0.014 \\
50  & 50 & 5 & 5 & \textbf{0.000} & 0.012 \\
50  & 50 & 8 & 2 & \textbf{0.000} & 0.006 \\
50  & 50 & 4 & 16 & \textbf{0.000} & 0.292 \\
50  & 50 & 10 & 10 & \textbf{0.000} & 0.255 \\
50  & 50 & 16 & 4 & $\mathbf{-0.001}$ & 0.145 \\
50  & 50 & 8 & 32 & \textbf{0.000} & 0.757 \\
50  & 50 & 20 & 20 & \textbf{0.000} & 0.771 \\
50  & 50 & 32 & 8 & \textbf{0.000} & 0.821 \\
100  & 100 & 4 & 16 & \textbf{0.000} & 0.043\\
100  & 100 & 10 & 10 & \textbf{0.000} & 0.037 \\
100  & 100 & 16 & 4 & \textbf{0.000} & 0.021 \\
100  & 100 & 8 & 32 & \textbf{0.000} & 0.762 \\
100  & 100 & 20 & 20 & \textbf{0.000} & 0.725 \\
100  & 100 & 32 & 8 & \textbf{0.000} & 0.527 \\
100  & 100 & 16 & 64 & \textbf{0.000} & 0.560 \\
100  & 100 & 40 & 40 & \textbf{0.000} & 0.584 \\
100  & 100 & 64 & 16 & \textbf{0.000} & 0.667 \\
200  & 200 & 8 & 32 & $\mathbf{-0.001}$ & 0.182 \\
200  & 200 & 20 & 20 & \textbf{0.000} & 0.157 \\
200  & 200 & 32 & 8 & \textbf{0.000} & 0.088 \\
200  & 200 & 16 & 64 & \textbf{0.000} & 0.807 \\
200  & 200 & 40 & 40 & \textbf{0.000} & 0.817 \\
200  & 200 & 64 & 16 & \textbf{0.000} & 0.850 \\
200  & 200 & 32 & 128 & \textbf{0.000} & 0.163 \\
200  & 200 & 80 & 80 & \textbf{0.000} & 0.210 \\
200  & 200 & 128 & 32 & \textbf{0.000} & 0.387 \\
50  & 100 & 2 & 8 & \textbf{0.000} & 0.013 \\
50  & 100 & 5 & 5 & \textbf{0.000} & 0.011 \\
50  & 100 & 8 & 2 & \textbf{0.000} & 0.005 \\
50  & 100 & 4 & 16 & \textbf{0.000} & 0.291 \\
50  & 100 & 10 & 10 & \textbf{0.001} & 0.244 \\
50  & 100 & 16 & 4 & \textbf{0.000} & 0.112 \\
50  & 100 & 8 & 32 & \textbf{0.000} & 0.757 \\
50  & 100 & 20 & 20 & \textbf{0.000} & 0.776 \\
50  & 100 & 32 & 8 & \textbf{0.000} & 0.840 \\
100  & 200 & 4 & 16 & \textbf{0.000} & 0.043 \\
100  & 200 & 10 & 10 & \textbf{0.000} & 0.035 \\
100  & 200 & 16 & 4 & \textbf{0.000} & 0.016 \\
100  & 200 & 8 & 32 & $\mathbf{-0.001}$ & 0.760 \\
100  & 200 & 20 & 20 & \textbf{0.000} & 0.711 \\
100  & 200 & 32 & 8 & \textbf{0.001} & 0.433 \\
100  & 200 & 16 & 64 & \textbf{0.000} & 0.560 \\
100  & 200 & 40 & 40 & \textbf{0.000} & 0.592 \\
100  & 200 & 64 & 16 & \textbf{0.000} & 0.711 \\
\bottomrule
\end{tabular}
\end{table}
\begin{table}[H]
\centering
\ContinuedFloat
\caption{(Continued)}
\setlength{\tabcolsep}{20pt}
\begin{tabular}{rrrrrr}
\toprule
$N_1$ & $N_2$ & $p_1$ & $p_2$ & $B_{\mathrm{prop}}$ & $B_{\mathrm{SYS}}$ \\
\midrule
200  & 400 & 8 & 32 & \textbf{0.001} & 0.182 \\
200  & 400 & 20 & 20 & \textbf{0.000} & 0.149 \\
200  & 400 & 32 & 8 & \textbf{0.000} & 0.065 \\
200  & 400 & 16 & 64 & \textbf{0.000} & 0.807 \\
200  & 400 & 40 & 40 & \textbf{0.000} & 0.821 \\
200  & 400 & 64 & 16 & \textbf{0.000} & 0.851 \\
200  & 400 & 32 & 128 & \textbf{0.000} & 0.164 \\
200  & 400 & 80 & 80 & \textbf{0.000} & 0.225 \\
200  & 400 & 128 & 32 & \textbf{0.000} & 0.451 \\
50  & 25 & 2 & 8 & \textbf{0.000} & 0.014 \\
50  & 25 & 5 & 5 & \textbf{0.000} & 0.012 \\
50  & 25 & 8 & 2 & \textbf{0.000} & 0.009 \\
50  & 25 & 4 & 16 & \textbf{0.000} & 0.292 \\
50  & 25 & 10 & 10 & \textbf{0.000} & 0.268 \\
50  & 25 & 16 & 4 & \textbf{0.000} & 0.188 \\
50  & 25 & 8 & 32 & \textbf{0.000} & 0.757 \\
50  & 25 & 20 & 20 & \textbf{0.000} & 0.766 \\
50  & 25 & 32 & 8 & \textbf{0.000} & 0.801 \\
100  & 50 & 4 & 16 & \textbf{0.000} & 0.043 \\
100  & 50 & 10 & 10 & \textbf{0.000} & 0.039 \\
100  & 50 & 16 & 4 & \textbf{0.000} & 0.027 \\
100  & 50 & 8 & 32 & \textbf{0.000} & 0.763 \\
100  & 50 & 20 & 20 & \textbf{0.000} & 0.738 \\
100  & 50 & 32 & 8 & \textbf{0.000} & 0.619 \\
100  & 50 & 16 & 64 & \textbf{0.000} & 0.559 \\
100  & 50 & 40 & 40 & \textbf{0.000} & 0.576 \\
100  & 50 & 64 & 16 & \textbf{0.000} & 0.640 \\
200  & 100 & 8 & 32 & $\mathbf{-0.001}$ & 0.182 \\
200  & 100 & 20 & 20 & \textbf{0.000} & 0.165 \\
200  & 100 & 32 & 8 & \textbf{0.000} & 0.114 \\
200  & 100 & 16 & 64 & \textbf{0.000} & 0.806 \\
200  & 100 & 40 & 40 & \textbf{0.000} & 0.814 \\
200  & 100 & 64 & 16 & \textbf{0.000} & 0.838 \\
200  & 100 & 32 & 128 & \textbf{0.000} & 0.162 \\
200  & 100 & 80 & 80 & \textbf{0.000} & 0.194 \\
200  & 100 & 128 & 32 & \textbf{0.000} & 0.317 \\
\bottomrule
\end{tabular}
\end{table}

\bibliographystyle{apalike} 
\bibliography{References.bib}

@book{muirhead1982aspects,
  title={Aspects of Multivariate Statistical Theory},
  author={Muirhead, Robb J},
  year={1982},
  publisher={John Wiley \& Sons, New York}
}

@book{Anderson2003aspects,
  author    = "Anderson, Theodore Wilbur",
  title     = "An introduction to multivariate statistical analysis",
  publisher = "Wiley-Interscience",
  year      = "2003",
  edition   = "3rd",
  series    = "Wiley series in probability and mathematical statistics",
  number    = "",
  URL       = "https://ci.nii.ac.jp/ncid/BA63301060"
}

@article{akita2010high,
title = {High-dimensional {E}dgeworth expansion of a test statistic on independence and its error bound},
journal = {Journal of Multivariate Analysis},
volume = {101},
number = {8},
pages = {1806-1813},
year = {2010},
issn = {0047-259X},
doi = {https://doi.org/10.1016/j.jmva.2010.03.014},
url = {https://www.sciencedirect.com/science/article/pii/S0047259X10000734},
author = {Tomoyuki Akita and Jinghua Jin and Hirofumi Wakaki}
}

@article{chang2010finite,
title = {Finite-sample inference with monotone incomplete multivariate normal data, {II}},
journal = {Journal of Multivariate Analysis},
volume = {101},
number = {3},
pages = {603-620},
year = {2010},
issn = {0047-259X},
doi = {https://doi.org/10.1016/j.jmva.2009.09.011},
url = {https://www.sciencedirect.com/science/article/pii/S0047259X09001766},
author = {Wan-Ying Chang and Donald St. P. Richards}
}

@article{wakaki2007error,
  title={Error bounds for high-dimensional {E}dgeworth expansions for some tests on covariance matrices},
  author={Wakaki, Hirofumi},
  journal={Hiroshima Statistical Research Group Technical Report,(07-04)},
  year={2007}
}

@article{wakaki2006edgeworth,
title = {Edgeworth expansion of {W}ilks’ lambda statistic},
journal = {Journal of Multivariate Analysis},
volume = {97},
number = {9},
pages = {1958-1964},
year = {2006},
note = {Special Issue dedicated to Prof. Fujikoshi},
issn = {0047-259X},
doi = {https://doi.org/10.1016/j.jmva.2005.10.010},
url = {https://www.sciencedirect.com/science/article/pii/S0047259X06000790},
author = {Hirofumi Wakaki}
}

@article{Sato2025sphericity,
	author = {Sato, Tetsuya and Yagi, Ayaka and Seo, Takashi},
	date = {2025/02/26},
	doi = {10.1007/s42519-025-00431-9},
	id = {Sato2025},
	isbn = {1559-8616},
	journal = {Journal of Statistical Theory and Practice},
	number = {2},
	pages = {16},
	title = {Sphericity Test on Variance-Covariance Matrix with Monotone Missing Data},
	url = {https://doi.org/10.1007/s42519-025-00431-9},
	volume = {19},
	year = {2025}}

@article{satomulti,
  title={An extension of sphericity test to the multi-sample problem with monotone incomplete data},
  author={Sato, Tetsuya and Yagi, Ayaka and Seo, Takashi},
  journal={Journal of Statistical Theory and Practice},
 volume={20},
  number={2},
  pages={54},
  year={2026},
  publisher={Springer}
}

@article{wakaki2018computable,
author = {Wakaki, H. and Fujikoshi, Ya.},
title = {Computable Error Bounds for High-Dimensional Approximations of an {LR} Statistic for Additional Information in Canonical Correlation Analysis},
journal = {Theory of Probability \& Its Applications},
volume = {62},
number = {1},
pages = {157-172},
year = {2018},
doi = {10.1137/S0040585X97T98854X},

URL = { 
    
        https://doi.org/10.1137/S0040585X97T98854X
    
    

},
eprint = { 
    
        https://doi.org/10.1137/S0040585X97T98854X
    
    

}
}

@article{kato2010asymptotic,
title = {High-dimensional asymptotic expansion of {LR} statistic for testing intraclass correlation structure and its error bound},
journal = {Journal of Multivariate Analysis},
volume = {101},
number = {1},
pages = {101-112},
year = {2010},
issn = {0047-259X},
doi = {https://doi.org/10.1016/j.jmva.2009.05.006},
url = {https://www.sciencedirect.com/science/article/pii/S0047259X09001110},
author = {Naohiro Kato and Takayuki Yamada and Yasunori Fujikoshi}
}

@book{fujikoshi2020non,
  title={Non-asymptotic analysis of approximations for multivariate statistics},
  author={Fujikoshi, Yasunori and Ulyanov, Vladimir V},
  year={2020},
  publisher={Springer}
}

@article{Anderson1985mle,
title = {Maximum-likelihood estimation of the parameters of a multivariate normal distribution},
journal = {Linear Algebra and its Applications},
volume = {70},
pages = {147-171},
year = {1985},
issn = {0024-3795},
doi = {https://doi.org/10.1016/0024-3795(85)90049-7},
url = {https://www.sciencedirect.com/science/article/pii/0024379585900497},
author = {T.W. Anderson and I. Olkin}
}

@article{box1949general,
  title={A general distribution theory for a class of likelihood criteria},
  author={Box, George EP},
  journal={Biometrika},
  volume={36},
  number={3/4},
  pages={317--346},
  year={1949},
  publisher={JSTOR}
}

@article{kanda1998some,
  title={Some basic properties of the {MLE}'s for a multivariate normal distribution with monotone missing data},
  author={Kanda, Takashi and Fujikoshi, Yasunori},
  journal={American Journal of Mathematical and Management Sciences},
  volume={18},
  number={1-2},
  pages={161--190},
  year={1998},
  publisher={Taylor \& Francis}
}

@article{mauchly1940significance,
author = {John W. Mauchly},
title = {{Significance test for sphericity of a normal $n$-variate distribution}},
volume = {11},
journal = {The Annals of Mathematical Statistics},
number = {2},
publisher = {Institute of Mathematical Statistics},
pages = {204-209},
year = {1940},
doi = {10.1214/aoms/1177731915},
URL = {https://doi.org/10.1214/aoms/1177731915}
}

@article{john1972distribution,
    author = {John, S.},
    title = {The distribution of a statistic used for testing sphericity of normal distributions},
    journal = {Biometrika},
    volume = {59},
    number = {1},
    pages = {169-173},
    year = {1972},
    month = {04},
    issn = {0006-3444},
    doi = {10.1093/biomet/59.1.169},
    url = {https://doi.org/10.1093/biomet/59.1.169},
    eprint = {https://academic.oup.com/biomet/article-pdf/59/1/169/738164/59-1-169.pdf},
}

@article{gleser1966note,
author = {Leon J. Gleser},
title = {{A note on the sphericity test}},
volume = {37},
journal = {The Annals of Mathematical Statistics},
number = {2},
publisher = {Institute of Mathematical Statistics},
pages = {464-467},
year = {1966},
doi = {10.1214/aoms/1177699529},
URL = {https://doi.org/10.1214/aoms/1177699529}
}

@article{ledoit2002some,
author = {Olivier Ledoit and Michael Wolf},
title = {{Some hypothesis tests for the covariance matrix when the dimension is large compared to the sample size}},
volume = {30},
journal = {The Annals of Statistics},
number = {4},
publisher = {Institute of Mathematical Statistics},
pages = {1081-1102},
year = {2002},
doi = {10.1214/aos/1031689018},
URL = {https://doi.org/10.1214/aos/1031689018}
}

@article{wang2013sphericity,
author = {Qinwen Wang and Jianfeng Yao},
title = {{On the sphericity test with large-dimensional observations}},
volume = {7},
journal = {Electronic Journal of Statistics},
publisher = {Institute of Mathematical Statistics and Bernoulli Society},
pages = {2164--2192},
year = {2013},
doi = {10.1214/13-EJS842},
URL = {https://doi.org/10.1214/13-EJS842}
}

@article{batsidis2006k-step,
title = {Statistical inference for location and scale of elliptically contoured models with monotone missing data},
journal = {Journal of Statistical Planning and Inference},
volume = {136},
number = {8},
pages = {2606-2629},
year = {2006},
issn = {0378-3758},
doi = {https://doi.org/10.1016/j.jspi.2004.10.021},
url = {https://www.sciencedirect.com/science/article/pii/S0378375804004252},
author = {A. Batsidis and K. Zografos}
}

@article{john1971,
 ISSN = {00063444, 14643510},
 URL = {http://www.jstor.org/stable/2334322},
 author = {S. John},
 journal = {Biometrika},
 number = {1},
 pages = {123--127},
 publisher = {[Oxford University Press, Biometrika Trust]},
 title = {Some Optimal Multivariates Tests},
 urldate = {2026-02-11},
 volume = {58},
 year = {1971}
}

@article{Srivastava2006,
  title = {Some tests criteria for the covariance matrix with fewer observations than the dimension},
  volume = {10},
  ISSN = {1406-2283},
  journal = {Acta et Commentationes Universitatis Tartuensis de Mathematica},
  publisher = {University of Tartu},
  author = {Srivastava,  Muni S.},
  year = {2006},
  month = Dec,
  pages = {77–93}
}

@article{Wilks1938,
author = {S. S. Wilks},
title = {{The large-sample distribution of the likelihood ratio for testing composite hypotheses}},
volume = {9},
journal = {The Annals of Mathematical Statistics},
number = {1},
publisher = {Institute of Mathematical Statistics},
pages = {60 -- 62},
year = {1938}
}

@article{Wilks1932,
 author = {S. S. Wilks},
 journal = {Biometrika},
 number = {3/4},
 pages = {471--494},
 publisher = {[Oxford University Press, Biometrika Trust]},
 title = {Certain Generalizations in the Analysis of Variance},
 urldate = {2026-07-26},
 volume = {24},
 year = {1932}
}

@article{Ulyanov2016,
  title = {Non-Asymptotic Results for {C}ornish–{F}isher Expansions*},
  volume = {218},
  ISSN = {1573-8795},
  number = {3},
  journal = {Journal of Mathematical Sciences},
  publisher = {Springer Science and Business Media LLC},
  author = {Ulyanov,  V. V. and Aoshima,  M. and Fujikoshi,  Y.},
  year = {2016},
  month = {Sept},
  pages = {363–368}
}
\end{document}